\pdfoutput=1
\documentclass[11pt]{amsart}
\usepackage{graphicx}
\usepackage{subfigure}
\usepackage{amscd}
\usepackage[square,numbers,sort]{natbib}
\usepackage[small,bf]{caption}
\usepackage{color}
\usepackage[colorlinks,linkcolor=black,bookmarksopen,bookmarksnumbered,
citecolor=black,urlcolor=black]{hyperref}

\allowdisplaybreaks
\numberwithin{equation}{section}

\graphicspath{{figs/}}

\definecolor{darkgrn}{rgb}{0, 0.75, 0}

\def\R{\mathbb{R}}

\DeclareMathOperator{\curl}{curl}
\def\div{\operatorname{div}}

\newcommand{\vol}[1]{\abs{#1}}

\def\d{\operatorname{d}}

\DeclareMathOperator{\hodge}{\ast} 
\DeclareMathOperator{\dual}{\star}

\DeclareMathOperator{\boundary}{\partial}

\newcommand{\abs}[1]{\lvert#1\rvert}

\newcommand{\eval}[2]{\langle #1,#2 \rangle}

\def\v{v}
\def\g{g}
\def\n{\hat{n}}
\def\M{M}
\def\D{D}

\title[Darcy flow using Discrete  Exterior Calculus]
{Numerical method for Darcy flow \\ derived using Discrete
  Exterior Calculus}
\author[A. N. Hirani]{Anil N. Hirani} 
\address{Correspondence to: Professor Anil N. Hirani, 
  Department of Computer Science, 
  University of Illinois at Urbana-Champaign, 
  201 N. Goodwin Ave., Urbana, IL 61801.} 
\email{\href{mailto:hirani@cs.illinois.edu}{hirani@cs.illinois.edu}}
\urladdr{\url{http://www.cs.illinois.edu/hirani}}
\author[K. B. Nakshatrala]{Kalyana B. Nakshatrala}
\address{Professor Kalyana Babu Nakshatrala, 
  Department of Mechanical Engineering, 
  Texas A\&M University, College Station, TX 77843.} 
\email{knakshatrala@tamu.edu}
\author[J. H. Chaudhry]{Jehanzeb H. Chaudhry}
\address{Jehanzeb H. Chaudhry, 
  Department of Computer Science, 
  University of Illinois at Urbana-Champaign, 
  201 N. Goodwin Ave., Urbana, IL 61801.}
  \email{jhameed2@illinois.edu}
\date{}

\begin{document}

\begin{abstract}
  We derive a numerical method for Darcy flow, hence also for
  Poisson's equation in mixed (first order) form, based on discrete
  exterior calculus (DEC). Exterior calculus is a generalization of
  vector calculus to smooth manifolds and DEC is one of its
  discretizations on simplicial complexes such as triangle and
  tetrahedral meshes. DEC is a coordinate invariant discretization, in
  that it does not depend on the embedding of the simplices or the
  whole mesh. We start by rewriting the governing equations of Darcy
  flow using the language of exterior calculus. This yields a
  formulation in terms of flux differential form and pressure. The
  numerical method is then derived by using the framework provided by
  DEC for discretizing differential forms and operators that act on
  forms. We also develop a discretization for spatially dependent
  Hodge star that varies with the permeability of the medium. This
  also allows us to address discontinuous permeability. The matrix
  representation for our discrete non-homogeneous Hodge star is
  diagonal, with positive diagonal entries. The resulting linear
  system of equations for flux and pressure are saddle type, with a
  diagonal matrix as the top left block. The performance of the
  proposed numerical method is illustrated on many standard test
  problems. These include patch tests in two and three dimensions,
  comparison with analytically known solution in two dimensions,
  layered medium with alternating permeability values, and a test with
  a change in permeability along the flow direction. We also show
  numerical evidence of convergence of the flux and the pressure. A
  convergence experiment is also included for Darcy flow on a surface.
  A short introduction to the relevant parts of smooth and discrete
  exterior calculus is included in this paper. We also include a
  discussion of the boundary condition in terms of exterior calculus.
\end{abstract}


\maketitle

\section{Introduction}
We have discretized the equations of Darcy flow and obtained a
numerical method on staggered mesh pairs with fluxes and pressures
being the primary variables. The numerical method was obtained by
using \emph{discrete exterior calculus} (DEC) \cite{Hirani2003,
  DeHiLeMa2005, DeKaTo2008}. Exterior calculus generalizes vector
calculus to higher dimensions and to smooth manifolds
\cite{AbMaRa1988} and DEC is one of its discretizations. This
discretization yields numerical methods for solving partial
differential equations (PDEs) on simplicial complexes, such as
triangle, tetrahedral or higher dimensional simplicial meshes. A
recent implementation of DEC is described in \cite{BeHi2011}. DEC is
related to many discretizations of exterior calculus that have been
popular in or are currently being pursued in numerical
analysis. Others include finite element exterior calculus
\cite{ArFaWi2006, ArFaWi2010}, support operator method, and mimetic
and compatible discretizations of PDEs \cite{HySh1997a, HySh1997b,
  HySh1999, BrLiSh2005, BoHy2006, BrLiShSi2007}, the covolume method
\cite{Nicolaides1992, NiTr2006}, and staggered cell methods
like~\cite{Mattiussi1997, Mattiussi2002}. See~\cite{ArBoLeNiSh2006}
for a collection of recent papers in these fields. Those parts of
exterior calculus and DEC that are relevant to this paper are
summarized in Section~\ref{sec:DEC}.

The equations of Darcy flow model the flow of a viscous incompressible
fluid in a porous medium. The equations consist of Darcy's law (which
expresses force balance), the continuity equation, and the boundary
condition. The first two form a very simple pair of equations, being
Poisson's equation in first order form. In this case Darcy flow can be
rewritten as Poisson's equation with pressure as the unknown. In
applications however, it is often the velocity field or flux that is
of primary interest. In some methods for the single variable second
order formulation, there is a loss of smoothness and accuracy in going
from pressure to velocity although there are reconstruction techniques
that can recover the accuracy.

However, one important metric of accuracy is the local mass balance
property, which is highly desirable feature for a numerical
formulation especially for applications in flow through porous
media. The single field formulation (which is based on the
second-order form) does not have the local mass balance property, and
hence is not accurate with respect to this metric. On the other hand,
mixed formulations are based on the mixed form of the governing
equations, and tend to perform better with respect to local mass
balance. The proposed DEC formulation is based on the mixed form of
the governing equations, and by construction we have local mass
balance. In first order form (i.e., when Darcy law and the continuity
equations are not combined into a single equation) Darcy flow
equations can be discretized using mixed finite element or volume
methods.

We now describe the primary distinctions between DEC and the above
mentioned related methods like finite element exterior calculus,
Support Operator Method, mimetic and compatible discretizations,
covolume method, and the staggered cell methods.  DEC was designed to
be coordinate invariant in order to mimic the coordinate invariance of
calculus on manifolds~\cite{Hirani2003, AbMaRa1988} which is a primary
language of modern physics~\cite{Frankel2004}. The specific embedding
of each simplex and of the whole mesh does not affect the resulting
discrete algebraic equations formulated via DEC. This is because the
Hodge star is defined in terms of lengths, areas, volumes etc. of
pieces of simplices which stay invariant with respect to simplex or
mesh embedding. One main distinction from mimetic discretizations and
support operator method is that those are designed for Euclidean
domains and are dependent on the embedding. As a result, DEC can be
used on simplicial approximations of non-flat manifolds as has been
shown by others~\cite{MuMcPaDuToKaMaDe2010,PaMuToKaMaDe2009} and by
the numerical experiments on a surface included in this paper and
in~\cite{HiKa2011}. There are methods such as those
in~\cite{DuGuJu2003a,DuJu2005,RiJuGu2008,JuDu2009} that are designed
to work on specific surfaces. However these are neither mimetic, nor
do they use a two variable formulation for diffusion like our method
does. Some other surface finite element methods are
in~~\cite{Dziuk1988,DeDz2007,Demlow2009} but again these are not
mimetic discretization and not formulated in mixed form. The recent
work in~\cite{HoSt2010} develops an abstract framework for finite
element exterior calculus for surfaces using a variational crimes
framework. However, no specific finite elements are developed and no
numerical experiments are shown in~\cite{HoSt2010}. The work of
Mattiussi~\cite{Mattiussi1997, Mattiussi2002} introduced algebraic
topology concepts to numerical methods. However, it is again dependent
on the embedding in Euclidean domains. The work that is closest to DEC
is finite element exterior calculus~\cite{ArFaWi2006,
  ArFaWi2010}. Indeed, by replacing the DEC Hodge star of our method
by a Whitney Hodge star~\cite{BeHi2011} one can obtain a finite
element exterior calculus method for Darcy flow. This has been
explored experimentally in~\cite{HiKa2011}. As Bochev and Hyman point
out in~\cite{BoHy2006}, all mimetic methods are related in that they
all are cochain based methods. The difference is in how (and if) the
Hodge star is discretized. In DEC this is done via circumcentric
duality in a coordinate or embedding invariant manner.

\begin{figure}[t!]
  \centering
  \includegraphics[scale=0.4, trim = 0.5in 2.5in 1in 2in, clip]
  {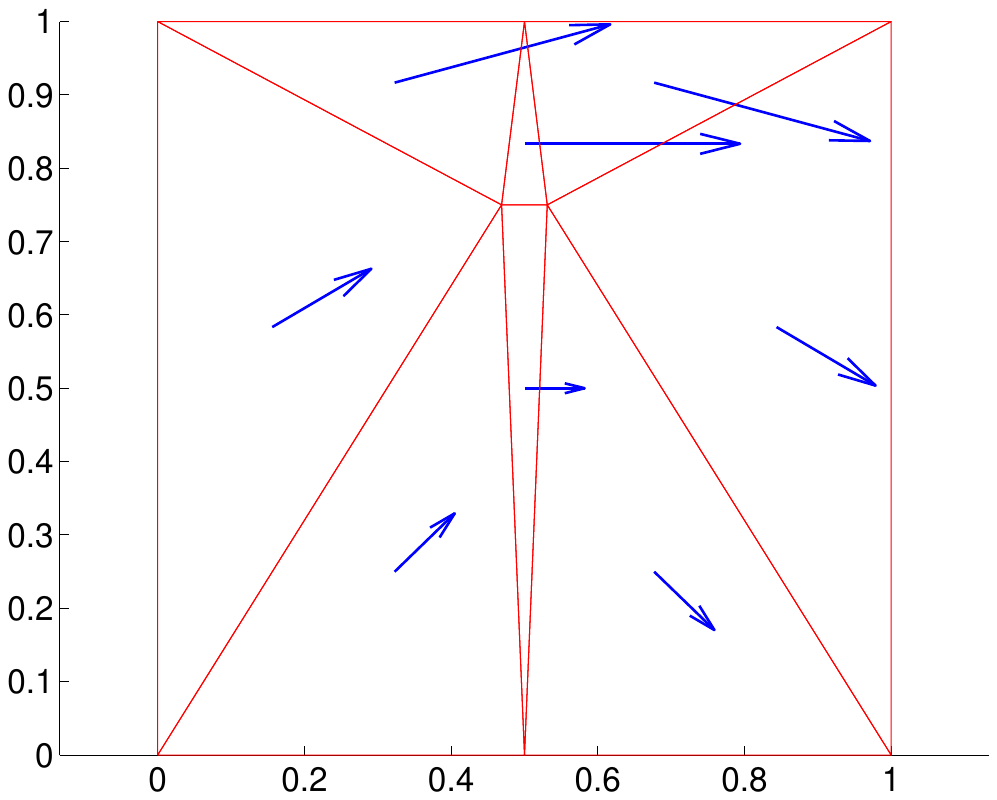}
  \caption{Mixed finite element method can be used to solve Darcy flow
    and it is well-known that care is needed in selecting the
    underlying finite element spaces. For example, equal order
    interpolation for both velocity and pressure is unstable. Here we
    show results for such a choice using piecewise linear finite
    elements. The correct solution is the constant vector field
    $(1,0)$. The fix for such unstable behavior is well-known in the
    finite elements literature \cite{BrFo1991, Braess2007}. In
    this paper we provide a related but different numerical
    method based on discrete exterior calculus.}
  \label{fig:sqr8nodal}
\end{figure}

In mixed methods, velocity (or flux) along with pressure are taken to
be the primary variables and this can yield more accurate results
compared to the pressure-only formulation. Mixed formulations require
careful choice of spaces for velocity and pressure since not all
combinations yield stable methods and the use of the LBB condition
\cite{BrFo1991, ArFaWi2010} is important here. For example,
the use of continuous piecewise linear representation for both
velocity and pressure results in an unstable method \cite{BrFo1991,
  Braess2007}. An example of such unstable behavior is shown in
Figure~\ref{fig:sqr8nodal}. Many fixes for such instability are
presented in the literature.  For example, one can use Raviart-Thomas
(RT) elements \cite{RaTh1977}, N\'ed\'elec elements
\cite{Nedelec1980}, Brezzi-Douglas-Marini (BDM) elements
\cite{BrDoMa1985}, or a variety of other finite element spaces
summarized in~\cite{Chen2005}.  Many of these spaces that have yielded
stable discretizations of Darcy flow and other problems have been
unified under the umbrella of finite element exterior calculus
\cite{ArFaWi2006}. Stabilized mixed finite elements methods have also
been developed for Darcy flow~\cite{MaHu2002,HuMaWa2006,BoDo2006,
  NaTuHjMa2006,NaMaHj2008}. Finite volume \cite{AcBeCo2003} and
covolume methods \cite{Nicolaides1992, ChKwVa1998,ChVa1999} are yet
another approach for solving the equations of Darcy flow.  A monotone,
locally conservative finite volume scheme for diffusion equation is
described in~\cite{LiShSvVa2007} and a mimetic finite difference
scheme is in~\cite{LiShSv2006}.

The primary variables in our method are area or volume flux (depending
on whether the problem is 2D or 3D), and pressure. The fluxes are
placed on edges in triangle meshes or triangles in tetrahedral meshes.
This leads to pressures being placed at circumcenters of top
dimensional simplices. There is numerical evidence that the
circumcenters need not be inside the simplices for the method to
work. See Figure~\ref{fig:sqr55cnstvel} for an example where the
circumcenters are not all inside the corresponding simplices. In all
our experiments involving high quality Delaunay meshes (which most
modern meshing software such as Triangle and tetgen automatically
provide~\cite{Shewchuk1996, Si2009}) the method has worked in the
presence of obtuse triangles. A sufficient condition for the method to
be applicable is that the meshes be well-centered (simplices should
contain their circumcenters). There has been considerable work by the
first author and collaborators on the creation of such
meshes. See~\cite{VaHiGuRa2010} and the references therein. Our method
is related to the methods for diffusion described in~\cite{PeSu2007}
and to the covolume method applied to Darcy flow
in~\cite{ChKwVa1998,ChVa1999} but our emphasis is on the relationships
between the numerical method and exterior calculus. In addition we
treat discontinuous permeability explicitly. The treatment of flux and
pressure in our method is similar to that for Navier-Stokes equations
in \cite{ElToKaScDe2007}, which is also based on DEC.

\subsection*{Our results: } Our discretization space is similar to the
lowest order Raviart-Thomas elements. However, the linear system
matrix resulting from our discretization is a saddle type matrix
\cite{BeGoLi2005} with a \emph{diagonal} matrix as the top left
block. Moreover our method is invariant with respect to coordinates or
embedding of the simplices since the Hodge star is computed from
quantities intrinsic to each simplex. Our method enforces mass balance
locally and globally and passes several standard numerical tests. In
2D we show that the method passes patch test involving constant
velocity and linear pressure. We also compare our numerical solution
with an analytical solution for a more general source term. We develop
a diagonal discretization of spatially dependent Hodge star that
depends on the permeability. That is, the matrix representation of the
Hodge star operator is diagonal. This is used for the case of layered
medium with alternating permeabilities, and for a discontinuous medium
where we use different pairs of permeability under constant velocity
conditions. We also show that our method passes patch tests in the 3D
case. We show numerical convergence of pressure and flux computed by
our method on a planar mesh. The numerically measured order of
convergence for flux is about 1.9 and that for pressure is about
1.04. A numerical convergence study is also done for flow on a surface
where the convergence orders are about 1.04 and 1 for flux and
pressure, respectively. All these numerical results are in
Section~\ref{sec:nmrclrslts} and the advantages of a DEC based
approach are discussed in Section~\ref{sec:conclusions}. A discussion
of boundary conditions in terms of exterior calculus is in
Sections~\ref{sec:gvrngeqtns} and~\ref{subsec:patch2d}.

\section{Review of Discrete Exterior Calculus (DEC)}
\label{sec:DEC}
In this section we briefly outline the relevant parts of smooth and
discrete exterior calculus. We discuss only the operators that are
relevant for Darcy flow. For more details on DEC
see~\cite{Hirani2003,DeHiLeMa2005,DeKaTo2008} and for details on
exterior calculus see~\cite{Arnold1989,AbMaRa1988}. In
Sections~\ref{sec:discretization} and \ref{sec:heterogeneous} we
describe our method for solving equations of Darcy flow so that it can
be implemented without knowledge of DEC or exterior calculus.  Thus, a
reader unfamiliar with some of the terms used in this section should
still be able to follow and implement the method. One useful
characteristic of exterior calculus is that all objects and operators
can be expressed in coordinate independent fashion. This aspect
however is harder to explain in a few paragraphs. Instead, we give
some examples using coordinates to describe the operators and objects
of exterior calculus.

\subsection{Smooth exterior calculus}
\label{subsec:smooth}
As mentioned earlier, exterior calculus generalizes vector calculus to
smooth manifolds \cite{AbMaRa1988,Arnold1989} and it consists of
operators on smooth general tensor fields defined on manifolds. A
tensor field evaluated at a point is a multilinear function on the
tangent space, mapping vector and covector arguments to $\mathbb{R}$
(the set of real numbers).  Other ranges besides $\mathbb{R}$ are
possible but not relevant here. Vector fields, symmetric tensor fields
such as metrics and antisymmetric tensor fields are all examples of
tensors. Antisymmetric tensors have been singled out in exterior
calculus and are called differential forms. A differential $k$-form
when evaluated at a point is an antisymmetric multilinear map on the
tangent space that takes $k$ vector arguments and produces a real
number. It is an object that can be integrated on a $k$-dimensional
space. In exterior calculus, it only makes sense to integrate
differential $k$-forms on a $k$-manifold.

Let $M$ be an $n$-dimensional orientable Riemannian manifold (a
manifold with an inner product on the tangent space at each point),
$TM$ the tangent bundle (disjoint union of the tangent spaces at all
points of $M$), $\mathfrak{X}(M)$ the space of smooth vector fields
and $\Omega^k(M)$ the space of differential $k$-forms on $M$.

Then \emph{exterior derivative} is a map $\d_k : \Omega^k(M)
\rightarrow \Omega^{k+1}(M)$ (sometimes written without the subscript)
that raises the degree of a form, and the \emph{wedge product} is a
map or binary operator $\wedge : \Omega^k(M) \times \Omega^l(M)
\rightarrow \Omega^{k+l}(M)$ that combines differential forms.  The
most important property of $\d$ is that $\d_{k+1} \circ \d_k =
0$. These two operators are enough to describe a basis for
differential forms on $M$. Taking $M = \mathbb{R}^3$ (the standard
three dimensional Euclidean space), with standard metric and
coordinates $x$, $y$ and $z$, a basis for $\Omega^1(\mathbb{R}^3)$,
the space of 1-forms is $(dx, dy, dz)$ and a basis for
$\Omega^2(\mathbb{R}^3)$ is $(dx \wedge dy, dx \wedge dz, dy \wedge
dz)$.  Let $f$ be a scalar valued function on $\mathbb{R}^3$ (i.e., a
$0$-form). Then its exterior derivative $\d f$ equals its differential
$df$ and is
\[
\d f = \frac{\partial f}{\partial x} dx +
                \frac{\partial f}{\partial y} dy +
                \frac{\partial f}{\partial z} dz\, .
\]
For a 1-form $\alpha = \alpha_1 dx + \alpha_2 dy + \alpha_3 dz$,
where $\alpha_i$ are scalar valued functions, its exterior derivative
is
\[
\d \alpha =
\left(
  \frac{\partial \alpha_2}{\partial x}-\frac{\partial \alpha_1}{\partial y}
\right) dx \wedge dy +
\left(
  \frac{\partial \alpha_3}{\partial x}-\frac{\partial \alpha_1}{\partial z}
\right) dx \wedge dz +
\left(
  \frac{\partial \alpha_3}{\partial y}-\frac{\partial \alpha_2}{\partial z}
\right) dy \wedge dz \, .
\]
The \emph{Hodge star} is an isomorphism, $\hodge : \Omega^k(M)
\rightarrow \Omega^{n-k}(M)$. For $\mathbb{R}^3$ with standard metric,
the Hodge star satisfies the following properties:
\begin{align*}
  &\hodge 1 = dx \wedge dy \wedge dz\,;\\
  &\hodge dx = dy \wedge dz,\quad
  \hodge dy = -dx \wedge dz,\quad
  \hodge dz = dx \wedge dy\,;\\
  &\hodge(dy \wedge dz) = dx, \quad
  \hodge(dx \wedge dz) = -dy, \quad
  \hodge(dx \wedge dy) = dz\,;\\
  &\hodge(dx \wedge dy \wedge dz) = 1 \, .
\intertext{For $\mathbb{R}^2$ the equivalent properties are}
  &\hodge 1 = dx \wedge dy\,;\\
  &\hodge dx = dy,\quad \hodge dy = -dx\,;\\
  &\hodge(dx \wedge dy) = 1 \, .
\end{align*}
An important property of Hodge star is that for a $k$-form $\alpha$,
\begin{equation}
  \label{eq:starstar}
  \hodge \hodge \alpha = (-1)^{k(n-k)}\alpha\, .
\end{equation}

Another operator relevant for Darcy flow is \emph{flat}, which is a
map $\flat : \mathfrak{X}(M) \rightarrow \Omega^1(M)$ that identifies
vector fields and 1-forms via the metric. Consider $\R^3$ with the
standard inner product, and standard orthonormal basis. Then for a
vector field $V$ with components $V_1$, $V_2$ and $V_3$, we have
$V^\flat = V_1 dx + V_2 dy + V_3 dz$. If the inner product is not the
standard one or the basis is not orthonormal then the relationship
between a vector field and its flat in coordinates is more
complicated. In $\mathbb{R}^3$, for a scalar function $f$ and vector
field $V$, some important relationships involving flat operator are:
\begin{equation}
(\nabla f)^\flat = \d f, \qquad
(\curl \, V)^\flat = \hodge \d V^\flat, \qquad
\div \, V = \hodge \d \hodge V^\flat\, .
\label{eq:vector_calculus_identities}
\end{equation}
Thus div, grad and curl can be defined in terms of exterior calculus
operators. Note that the operators $\d$ and $\wedge$ are metric
independent and so they can be defined on a manifold without having to
define a Riemannian metric. On the other hand, the operators $\hodge$
and $\flat$ do require a metric for their definition.

\subsection{Primal and dual mesh}
\label{subsec:prmldlmsh}
Discretizing exterior calculus involves deciding what should replace
the smooth manifolds, differential forms and other tensor fields and
operators that act on these. Recall that a \emph{simplicial complex}
$K$ in $\mathbb{R}^N$ is a collection of simplices in $\mathbb{R}^N$
such that every face of a simplex of $K$ is in $K$ and such that the
intersection of any 2 simplices of $K$ is a face of each of them. By
\emph{face} one means the simplex itself or one of its lower
dimensional subsimplices. For example, a triangle has as faces,
itself, the three edges, and the three vertices. The dimension $n$ of
the complex is the highest dimension of its simplices. Thus $n \le N$
where $N$ is the dimension of the embedding space. In DEC, the
oriented Riemannian manifolds $M$ of smooth exterior calculus is
replaced by the \emph{underlying space} $\abs{K}$ of an \emph{oriented
  manifold simplicial complex} $K$ as described in
\cite{Hirani2003}. Briefly, these are simplicial complexes in which
the neighborhood of every interior point is homeomorphic to (``looks
like'') an open subset of $\mathbb{R}^n$ and in which each simplex of
dimension $k$, for $0 \le k \le n-1$, is a face of some
$n$-simplex. Thus a triangle with an edge sticking out from one vertex
would not be admissible and neither would a triangle mesh surface with
a fin like triangle sticking out from an edge. 

In addition all the $n$-simplices must have the same orientation.  The
\emph{orientation} of a top dimensional simplex can be defined in the
same way that it is defined for an affine space, by specifying an
ordered tuple of basis vectors in the simplex as having one
orientation. For example one can take one vertex to be an origin and
the $n$ vectors to the remaining $n$ vertices as the tuple of
vectors. Equivalently, specifying one ordering of the vertices as one
orientation and any odd permutation of the ordering to be the other
orientation is also possible. For a manifold complex like a piecewise
triangular surface embedded in $\R^3$ the notion of triangles having
the same orientation is particularly simple. Each shared edge should be
traversed in opposite directions when following the orientation of one
triangle or the other. See \cite{Hirani2003} or~\cite{Munkres1984} for
details. An oriented manifold simplicial complex will be also called a
\emph{primal mesh}.

\begin{figure}[t]
  \centering
  \includegraphics[width=4in]{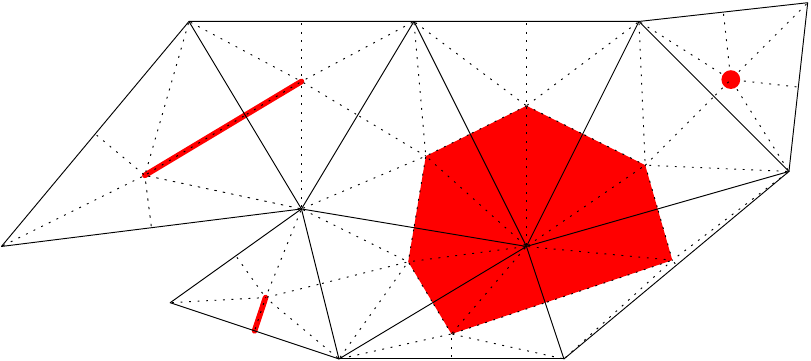}  
  \caption{A simplicial complex is subdivided using circumcentric
    subdivision and the dual cells are constructed from the
    subdivision. The new edges introduced by the subdivision are shown
    dotted. The dual cells shown are colored red. The dual cells are
    not explicitly stored, only their length, area, or volume is
    needed. Figure taken from~\cite{Hirani2003}.}
  \label{fig:primalDualComplex2d}
\end{figure}

In addition to the primal mesh $K$ a staggered cell complex $\dual K$
associated with $K$ and referred to as the \emph{dual mesh} also plays
a role in DEC. An example of a primal mesh, with some pieces of the
dual mesh highlighted is shown in
Figure~\ref{fig:primalDualComplex2d}. Usually the dual mesh is not
explicitly stored. What are needed instead are lengths, areas or
volumes of pieces of the dual mesh, the specific needed quantities
depending on the application. The primal mesh is a simplicial complex,
i.e., a triangle or tetrahedral (or higher dimensional) mesh such as
are used in finite volume or finite element methods. A sufficient
condition is that the mesh be well-centered~\cite{VaHiGuRa2010}. In
practice the class of meshes on which DEC applies is somewhere between
well-centered and Delaunay meshes, and a precise characterization is
not currently known.  In practice, good quality Delaunay meshes
produced by most modern meshing programs (such as Triangle
\cite{Shewchuk1996} and tetgen \cite{Si2009}) appear to suffice.

In what follows, $\sigma^k$ will be a \emph{primal} $k$-simplex, a
$k$-dimensional simplex in the primal mesh.  The corresponding
\emph{dual} $(n-k)$-cell, an $(n-k)$-dimensional cell in the dual mesh
will be denoted by $\dual \sigma^k$. We will use
$\operatorname{cc}(\sigma^k)$ to mean the circumcenter of $\sigma^k$,
the unique point equidistant from all vertices of $\sigma^k$. The
notation $\sigma \prec \tau$ will mean that simplex $\sigma$ is a face
of simplex $\tau$ and $\tau \succ \sigma$ will mean that $\tau$
contains $\sigma$ as one of its faces.

To find the dual cell $\dual \sigma^k$ of the primal simplex
$\sigma^k$ proceed as follows. Start from the circumcenter
$\operatorname{cc}(\sigma^k)$. Traverse in straight lines, one by one,
to the circumcenters $\operatorname{cc}(\sigma^{k+1})$ of all
$\sigma^{k+1} \succ \sigma^k$ and from those to the next dimension and
so on all the way to the top dimension. Each path of these traversals
yields a simplex of dimension $n-k$.  The union of all such simplices
is the dual cell $\dual \sigma^k$ of simplex $\sigma^k$.  For example,
in Figure~\ref{fig:primalDualComplex2d}, the dual of an internal edge
is obtained by starting from its middle and traversing to the
circumcenters of the two triangles containing it. The resulting two
edges together form the dual of the edge. Note that it need not be a
straight line if the two adjacent triangles do not lie in the same
affine space. See~\cite{Hirani2003} for more examples of primal-dual
pairs.

The primal and dual meshes of DEC are oriented. The primal mesh is
oriented consistently at the top level. For example, either all the
triangles in a triangle mesh must be oriented clockwise, or all of
them must be oriented counter-clockwise. Similarly, all the tetrahedra
in a tetrahedral mesh must be right-handed, or all must be
left-handed. The lower dimensional simplices can be oriented
arbitrarily, for example, using the dictionary order of the vertex
numbers. The orientation of the dual cells is implied by the
orientation of the corresponding primal simplices and of the top level
primal simplices.  For details see~\cite{Hirani2003}.  For triangle
meshes, a vector along the primal edge followed by one along the dual
edge should define the same orientation as that of the triangles. For
tetrahedral meshes the orientation of a face followed by a vector
along the dual edge should form the same handedness as that of the
tetrahedron.

\subsection{Chains and boundary operator}
\label{subsec:chains}
Let $K$ be a finite simplicial complex. Recall from algebraic topology
\cite{Munkres1984} that a $k$-\emph{chain} on $K$ is a function from
the set of oriented $k$-simplices of $K$ to the set of integers
$\mathbb{Z}$. A $k$-chain $c$ has the property that $c(\sigma) =
-c(-\sigma)$, where $-\sigma$ is $\sigma$ with the opposite
orientation. Chains are added by adding their values. The space of
$k$-chains is a group and is denoted $C_k(K)$. The group structure
will not be important to us except for the fact that it will allow us
to talk about homomorphisms~--~maps between groups that preserve the
group structure. For a $k$-simplex $\sigma^k$ we will use $\sigma^k$
or $\sigma$ to denote the simplex as well as the chain that takes the
value 1 on the simplex and 0 on all other $k$-simplices in $K$. Such a
chain is called an \emph{elementary $k$-chain}.

The \emph{boundary operator} $\boundary_k : C_k(K) \rightarrow
C_{k-1}(K)$ is defined as a homomorphism (that is, $\boundary_k(a+b)
= \boundary_k a + \boundary_k b$) by defining it on oriented simplices
$[v_0, \ldots v_k]$. It is defined by
\[
\boundary_k [v_0 \ldots v_k] = 
\sum_{i=0}^k(-1)^i[v_0,\ldots,\widehat{v_i},\ldots,v_k]
\]
the hat indicating the missing vertex. For example if $[v_0,v_1,v_2]$
is a triangle in $\mathbb{R}^2$, oriented counterclockwise, then the
boundary of the corresponding elementary 2-chain is the sum of the 3
elementary 1-chains $[v_1,v_2]$, $-[v_0,v_2](=[v_2,v_0])$ and
$[v_0,v_1]$.

\subsection{Cochains as discrete forms}
\label{subsec:dscrtfrms}
As is usual in most discretizations of exterior calculus, in DEC
discrete differential $k$-forms are defined to be elements of
$\operatorname{Hom}\bigl(C_k(K), \mathbb{R}\bigr)$. This is the space
of real-valued homomorphisms on the space of $k$-chains. This space is
called the space of $k$-\emph{cochains} or \emph{discrete differential
  $k$-forms}. Thus given $k$-chains $a$ and $b$ and a $k$-cochain
$\alpha$, we have $\alpha(a+b) = \alpha(a) + \alpha(b)$ in which each
term is a real number. The space of primal $k$-cochains on a
simplicial complex $K$ will be denoted by $C^k(K;\mathbb{R})$ and the
dual $k$-cochains on the dual cell complex $\dual K$ by $C^k(\dual
K;\mathbb{R})$. We will shorten these notations to $C^k(K)$ and
$C^k(\dual K)$. On the underlying space $\abs{K}$ of the simplicial
complex $K$, we will be working with piecewise smooth differential
forms. More precisely DEC starts with square integrable forms
\cite{Dodziuk1976} on $\abs{K}$, denoted as $L^2\Omega^k(\abs{K})$ and
so that is the notation we will use for piecewise smooth
forms. Discrete forms are created from (piecewise) smooth forms with
the \emph{de Rham} map $R : L^2\Omega^k(\abs{K}) \rightarrow C^k(K)$
or $R : L^2\Omega^k(\abs{K}) \rightarrow C^k(\dual K)$ depending on
the context. See \cite{Dodziuk1976} for details on the de Rham
map. For a form $\alpha$, we will denote the evaluation of the cochain
$R(\alpha)$ on a chain $c$ as $\eval{R(\alpha)}{c}$ and define it as
$\int_c \alpha$. Thus given a smooth $k$-form $\alpha$ the de Rham map
converts it into the $k$-cochain $\int_{\text{\textvisiblespace}}
\alpha$ with the slot for integration domain left empty. This cochain
$\int_{\text{\textvisiblespace}} \alpha$ is ready to be applied to a
$k$-chain $c$ to produce a number $\int_c \alpha$. Note that we are
implicitly assuming that the smooth quantities are defined on the
simplicial mesh that is the discretization of the smooth manifold. For
planar and spatial domains we consider in the examples in this paper,
this condition is trivially true. For more general domains, like
surfaces embedded in $\mathbb{R}^3$ this restriction can be removed,
but work on such generalizations, especially how it affects numerical
solutions of PDEs is still in early stages. For some related ideas see
\cite{HiPoWa2006, Wardetzky2008, DeDz2007, Demlow2009, HoSt2010}. For
the Darcy flow problem, the pressure is a dual 0-cochain and the flux
will be a primal 1-cochain for 2D problems and a 2-cochain for 3D
problems.

The cochain spaces $C^k(K)$ and $C^k(\dual K)$ defined above are
groups, but in addition they are also vector spaces. If there are
$N_k$ simplices of dimension $k$ in $K$ then the vector space
dimension $\dim\bigl(C^k(K)\bigr)$ is $N_k$. Similarly, if there are
$N_k$ cells of dimension $k$ in $\dual K$ then $\dim\bigl(C^k(\dual
K)\bigr)$ is $N_k$. Note that $\dim\bigl(C^k(K)\bigr) =
\dim\bigl(C^{n-k}(\dual K)\bigr)$ since $k$-simplices of $K$ are in
one-to-one correspondence with $(n-k)$-cells of $\dual K$. To define a
basis for $C^k(K)$ as a vector space, the $k$-simplices are first
ordered in some way. For example in PyDEC \cite{BeHi2011}, the
$k$-simplices are ordered in dictionary order based on the vertex
names of the simplex. In a triangle $[v_0, v_1, v_2]$, the dictionary
order for edges is $[v_0,v_1]$, $[v_0,v_2]$ and $[v_1,v_2]$. We would
typically refer to these edges as $\sigma^1_0$, $\sigma^1_1$ and
$\sigma^1_2$, respectively. If $\sigma_0, \ldots, \sigma_{N_k}$ is an
ordered list of the $k$-simplices of $K$ then we can define a basis
$\bigl(\sigma_0^\ast,\ldots,\sigma^\ast_{N_k}\bigr)$ for $C^k(K)$
where $\eval{\sigma_i^\ast}{\sigma_j} = \delta_{ij}$, the Kronecker
delta. That is, $\sigma_i^\ast$ is the $k$-cochain that is 1 on the
elementary $k$-chain of $\sigma_i$ and 0 on the other elementary
$k$-chains in $K$. In computations we represent elements of $C^k(K)$
and $C^k(\dual K)$ as vectors of appropriate dimensions in these
bases.

\subsection{Discrete exterior derivative and Hodge star}
\label{subsec:dscrtops}
We now give definitions of \emph{discrete exterior derivative} and
\emph{discrete Hodge star}. These are given here without explanation
as to why these are good choices for the discrete operators. Such an
explanation can be found in~\cite{Hirani2003,DeHiLeMa2005}.
The discrete exterior derivative on the primal cochains will also be
denoted as $\d$ (or $\d_k$ if the degree of the source space is to be
specified) and is defined using the boundary operator in such a way
that Stokes theorem is true by definition. For the exterior derivative
on the dual cochains we will use the notation $\d^\ast$ (or
$\d_k^\ast$). For a $k$-cochain $\alpha^k$ and $(p+1)$-chain $c^{k+1}$, define
\begin{equation}
\label{eq:discrete_d}
\eval{\d_k \alpha^k}{c^{k+1}} := \eval{\alpha^k}{\boundary_{k+1} c^{k+1}}\, .
\end{equation}

In the basis for $C^k(K)$ described in Section~\ref{subsec:dscrtfrms},
the matrix representation of $\d_k$, denoted $\D_k$, is an $N_{k+1}$ by
$N_k$ matrix with entries $0$, $1$ or $-1$. Analogous bases for
$C^k(\dual K)$ using the dual cells yield matrix form of the dual
discrete exterior derivative operator. The matrices of the primal and
dual exterior derivative are related. In fact the matrix form of
$\d_k^\ast$ is $\pm \D_{n-k-1}^T$ where the sign depends on $n$
and~$k$.  For Darcy flow the only relevant pairs of primal and dual
exterior derivatives are $\d_0^\ast$ and $\d_{n-1}$ and the matrix
form for $\d_0^\ast$ is $(-1)^n\D_{n-1}^T$.

The discrete Hodge star can be thought of as an operator that
``transfers information'' between the primal and dual meshes. Given a
$k$-simplex $\sigma$ and a primal $k$-cochain $\alpha$, the discrete
Hodge star of $\alpha$ (denoted $\hodge \alpha$) is a dual discrete
$(n-k)$-cochain defined by its value on the dual cell $\dual \sigma$
by
\begin{equation}
\label{eq:discrete_hodge}
\frac{1}{\vol{\dual\sigma}}\eval{\hodge \alpha}{\dual \sigma}
:= \frac{1}{\vol{\sigma}}\eval{\alpha}{\sigma}\, .
\end{equation}
Here $\vol{\sigma}$ is the measure of $\sigma$ and $\vol{\dual\sigma}$
is the measure of the circumcentric dual cell corresponding to
$\sigma$, with measure of a 0-dimensional object being~1. See
\cite{Hirani2003, DeHiLeMa2005} for details. We will sometimes use
$\hodge_k$ to denote the Hodge star with $C^k(K)$ as its domain. Using
the bases for $C^k(K)$ and $C^{n-k}(\dual K)$ mentioned above, the
matrix $\M_k$ for $\hodge_k$ is a diagonal $N_k$ by $N_k$ matrix.  It
is helpful to use $\hodge_k^{-1}$ to denote the inverse map although
the inverse notation is usually not used for smooth Hodge star. To
simulate the Hodge star property~\eqref{eq:starstar} on the discrete
side we will always use \begin{equation}
  \label{eq:dscrtstrstr}
  (-1)^{k(n-k)} \hodge_k^{-1} \, ,
\end{equation}
whenever the inverse discrete Hodge star map is used.

The various primal and dual cochain complexes in 2 dimensions are
related via the discrete exterior derivative and discrete Hodge star
operators in the manner shown in the following diagram:
\begin{equation}
  \label{eq:cochncmplx2d}
  \begin{CD}
    C^0(K) @>\d_0>> C^1(K) @>\d_1>> C^2(K) \\
    @VV\hodge_0 V    @VV\hodge_1 V    @VV\hodge_2 V \\
    C^2(\dual K) @<\d_1^\ast< < C^1(\dual K) @<\d_0^\ast< < 
    C^0(\dual K)
  \end{CD}
\end{equation}
The corresponding diagram for 3 dimensions is:
\begin{equation}
  \label{eq:cochncmplx3d}
  \begin{CD}
    C^0(K) @>\d_0>> C^1(K) @>\d_1>> C^2(K) @>\d_2>> C^3(K)\\
    @VV\hodge_0 V    @VV\hodge_1 V    @VV\hodge_2 V @VV\hodge_3 V \\
    C^3(\dual K) @<\d_2^\ast< < C^2(\dual K) @<\d_1^\ast< < 
    C^1(\dual K) @<\d_0^\ast< < C^0(\dual K)
  \end{CD}
\end{equation}

\subsection{Interpolation of cochains}
\label{subsec:intrp}
To go from cochains to piecewise smooth forms on a simplicial complex
$K$ a well-known map is the \emph{Whitney map}
\cite{Whitney1957,Bossavit1988a,Bossavit1988b} denoted $W$. This map
can be thought of as a way to interpolate numbers defined on edges,
triangles and tetrahedra. Thus it goes in direction opposite of the
one for the de Rham map which is a discretization map. For a complex
$K$ consisting of a triangle $[v_0, v_1, v_2]$ and its faces, the
Whitney map for 1-cochains is defined by extending by linearity, the
following to all of $C^1(K)$
\begin{equation}
  \label{eq:whitney_1form}
  W\bigl([v_i,v_j]^\ast\bigr) := \mu_i \: \d \mu_j - \mu_j 
  \: \d \mu_i \,  , 
\end{equation}
for $i \ne j$ and $i,j \in \{0,1,2\}$. Here $\mu_{i}$ is the
barycentric basis function corresponding to $v_i$, i.e. the affine
function that is 1 on $v_i$ and 0 at other vertices. Recall that the
1-cochain $[v_i,v_j]^\ast$ is 1 on the elementary 1-chain of the edge
$[v_i, v_j]$ and 0 at all other elementary 1-chains. For a tetrahedron
$[v_0, v_1, v_2, v_3]$ there are Whitney maps as above for
interpolating the edge values and in addition there are Whitney maps
for interpolating the triangle values, and these are defined by
extension from
\begin{equation}
  \label{eq:whitney_2form}
  W\bigl([v_i,v_j,v_k]^\ast\bigr) := 2 \bigl(\mu_i \: 
  \d \mu_j \wedge \d \mu_k - 
  \mu_j \: \d \mu_i \wedge \d \mu_k + \mu_k \: \d \mu_i 
  \wedge \d \mu_j\bigr) \, .
\end{equation}

The piecewise smooth forms (smooth in each simplex) constructed using
the Whitney map are also known as a \emph{Whitney forms}.  Whitney forms can
be used to build a low order finite element exterior calculus
\cite{Bossavit1988, Bossavit1988a} and in computational
electromagnetism the Whitney 1-form and 2-form are also known as edge
and face elements respectively \cite{Bossavit1998}. Finite element
exterior calculus has now been generalized to general polynomial
differential forms \cite{ArFaWi2006}. We use the Whitney maps only for
interpolating the differential forms so we can plot the corresponding
vector field for visualization. The vector field corresponding to the
Whitney 1-form in equation~\eqref{eq:whitney_1form} is obtained by
applying a sharp operator to get
\[
\mu_i \: \nabla \mu_j - \mu_j \: \nabla \mu_i \, .
\]
The vector field corresponding to the Whitey 2-form in 
equation~\eqref{eq:whitney_2form} is 
\[
2 \bigl(\mu_i \: \nabla \mu_j \times \nabla \mu_k - 
  \mu_j \: \nabla \mu_i \times \nabla \mu_k + 
  \mu_k \: \nabla \mu_i \times \nabla \mu_j\bigr) \, .
\]

\section{Governing Equations}
\label{sec:gvrngeqtns}
We first present the governing equations of Darcy flow in the standard
vector calculus notation, and then rewrite them using differential
forms and vector fields (that is, in exterior calculus notation). This
latter form is then discretized on a simplicial complex and its dual,
which yields a numerical method for Darcy flow.

Let $M \subset \mathbb{R}^n$ be a bounded open domain,
$\bar{M}$ its closure and $\partial M := \bar{M}
\backslash M$ its boundary, which is assumed to be piecewise
smooth. In this paper $n$ (which represents spatial dimensions) can be
2 or 3. Let $\v : M \rightarrow \mathbb{R}^n$ be the Darcy
velocity \cite{MaHu2002} (units  $\text{m}^2/(\text{m}\,\text{s}) =
\text{m}/\text{s}$ for $n=2$ or $\text{m}^3/(\text{m}^2\,\text{s}) =
\text{m}/\text{s}$ for $n=3$) and let $p : M \rightarrow
\mathbb{R}$ be the pressure. The governing equations of Darcy flow can
be written as \begin{alignat}{2}
  \v + \frac{\kappa}{\mu} \nabla p &= \frac{\kappa}{\mu} \rho \g 
  &&\quad \text{in } M \label{eq:vfdrcy}\, ,\\
  \operatorname{div}\v &= \phi 
  &&\quad \text{in } M \label{eq:vfcnt}\, ,\\
  \v \cdot \n &= \psi 
  &&\quad \text{on } \partial M\, ,
  \label{eq:vfbc}
\end{alignat}
where $\kappa > 0$ is the coefficient of permeability of the medium (units
$\text{m}^2$ for $n=3$), $\mu > 0$ is the coefficient of (dynamic)
viscosity of the fluid (units $\text{kg}/(\text{m}\,\text{s})$), $\rho
> 0$ is the density of the fluid, $\g$ is the acceleration due to
externally applied body force (i.e., $\rho \g$ is the body force
density), $\phi:M \rightarrow \mathbb{R}$ is the prescribed
divergence of velocity, $\psi:\partial M \rightarrow \mathbb{R}$
is the prescribed normal component of the velocity across the
boundary, and $\n$ is the unit outward normal vector to $\partial
M$. For consistency $\int_M \phi \, dM = \int_{\partial
  M} \psi \, d\Gamma$ where $d\Gamma$ is the measure on $\partial M$.

Equation \eqref{eq:vfdrcy} is Darcy's law, equation \eqref{eq:vfcnt}
is the continuity equation and equation \eqref{eq:vfbc} is the
boundary condition.  In the above equations, permeability $\kappa$ is
assumed to be a scalar constant. In Section~\ref{sec:heterogeneous} we
will relax this constraint and allow $\kappa$ to be a scalar valued
function of space. The further generalization needed for modeling
anisotropic permeability requires $\kappa$ to be a tensor, which is
not addressed in this paper.  To simplify the treatment, in the rest
of this paper we will assume that there is no external force acting on
the system, i.e., $\g = 0$.

The first step in the DEC formulation is to rewrite the governing
equations \eqref{eq:vfdrcy} -- \eqref{eq:vfbc} in exterior calculus
notation. As above, we first assume that the permeability $\kappa$ is a
scalar constant. We first apply the flat operator to both sides of
equation~\eqref{eq:vfdrcy}, use
equation~\eqref{eq:vector_calculus_identities} for divergence, and
then apply Hodge star to both sides of equation~\eqref{eq:vfcnt} to
obtain (assuming $\g = 0$)
\begin{alignat}{2}
    \v^{\flat} + \frac{\kappa}{\mu} \d_0 p &= 0
    &&\quad \text{in } M \label{eq:vflatdrcy}\, ,\\
    \d_{n-1}(\hodge\v^{\flat}) &= \phi \omega 
    &&\quad \text{in } M \label{eq:vflatcnt}\, ,\\
    \hodge\v^{\flat} &= \psi \gamma 
    &&\quad \text{on } \partial M\, .
    \label{eq:vflatbc}
\end{alignat}
Here $\omega = \hodge 1$ is a volume $n$-form on $M$ and $\gamma$
is the volume $(n-1)$-form on the boundary $\partial M$ and it is
defined by requiring
\begin{equation}
  \label{eq:gamma}
  \gamma(X_1,\ldots,X_{n-1}) = \omega(\n, X_1,\ldots, X_{n-1})\, ,
\end{equation}
for all vector fields $X_1,\ldots,X_{n-1}$ on the boundary $\partial
M$. Note that in going from equation~\eqref{eq:vfbc}
to~\eqref{eq:vflatbc} we have used the fact that $(\v \cdot \n)\gamma
= \hodge \v^\flat$. For an explanation of why this is true see \cite[page
506]{AbMaRa1988}. The definition~\eqref{eq:gamma} of
$\gamma$ has implications on how the orientations affect the sign of
the quantity $\psi \gamma$ and this is explained using a concrete
example in Section~\ref{sec:nmrclrslts}. The other quantities in the
equations above are as in~\eqref{eq:vfdrcy} -- \eqref{eq:vfbc}.  Note
that $\phi$ and $\psi$ must satisfy 
\[
\int_M \phi \omega =
\int_{\partial M} \psi \gamma\, ,
\]
by Stokes' theorem, which is analogous to the consistency condition
stated earlier.

Next we define a differential form which is the volumetric (or volume)
flux in 3D (units $\text{m}^3/(\text{m}^2 \, \text{s}) =
\text{m}/\text{s}$) or area flux in 2D (units $\text{m}^2/(\text{m} \,
\text{s}) = \text{m}/\text{s}$). This quantity will be denoted as $f$,
which is an $(n - 1)$-form, defined by
\[ 
f := \hodge (\v^{\flat}) \, . 
\] 
This is appropriate because as mentioned above, $(\v \cdot \n)\gamma =
\hodge \v^\flat$. Applying Hodge star to both sides of equation
\eqref{eq:vflatdrcy} and replacing $\hodge\v^{\flat}$ by $f$
everywhere we get the governing equations in terms of the flux and
pressure, which can be written as \begin{alignat}{2}
  f + \frac{\kappa}{\mu} (\hodge \d_0 p) &= 0
  &&\quad \text{in } M \label{eq:dfdrcy}\, ,\\
  \d_{n-1} f  &= \phi\omega 
  &&\quad \text{in } M \label{eq:dfcnt}\, ,\\
  f &= \psi\gamma 
  &&\quad \text{on } \partial M \label{eq:dfbc}\, .
\end{alignat}
Given $\kappa$, $\mu$, $\phi$, $\psi$ and the boundary condition
\eqref{eq:dfbc}, the problem statement is to solve equations
\eqref{eq:dfdrcy} and~\eqref{eq:dfcnt} for the flux $f$ and pressure
$p$.  Equations~\eqref{eq:dfdrcy}--\eqref{eq:dfbc} are the ones that
we will discretize first in Section~\ref{sec:discretization} using the
discrete operators defined in Section~\ref{subsec:dscrtfrms}.  An
equivalent form for equation~\eqref{eq:dfdrcy} obtained by applying
Hodge star to both sides of~\eqref{eq:dfdrcy} is
\begin{equation}
  \hodge f + (-1)^{n-1}\,\frac{\kappa}{\mu} \d_0 p = 0 
  \quad \text{in} \; M \label{eq:starfdrcy}\, ,
\end{equation}
and we will also discretize this equation to get an alternative
formulation. Here the $(-1)^{n-1}$ sign has come from the double
application of Hodge star using equation~\eqref{eq:starstar}.

\section{Discretization of Equations}
\label{sec:discretization}
Let $K$ be a simplicial complex that approximates $M$ and $\dual
K$ the circumcentric dual of $K$ as defined in
Section~\ref{subsec:prmldlmsh}. Let $L$ be the approximation of the
boundary $\partial M$ so that $L$ consists of the
$(n-1)$-dimensional boundary faces of $K$. The differential forms $f$,
$k$, $\phi\omega$ and $\psi \gamma$ in
equations~\eqref{eq:dfdrcy}-~\eqref{eq:dfbc} are discretized as
cochains and the operators $\d$ and $\hodge$ are replaced by their
discrete counterparts described in Section~\ref{subsec:dscrtops}.

An important point to note when discretizing is the appropriate
placement of the cochains. In particular we will place the discrete
flux $f$ on the $(n-1)$-dimensional \emph{primal} simplices -- edges
in triangle mesh and triangles in tetrahedral mesh. Thus $f \in
C^{n-1}(K)$. From equation~\eqref{eq:dfdrcy} this implies that the
discrete version of $\hodge \d_0 p$ must also be placed on these primal
simplices since $k/\mu$ is a scalar here. Thus $\hodge \d_0 p \in
C^{n-1}(K)$ from which it follows that $\d_0 p$ is a dual cochain and
$\d_0 p \in C^{n-(n-1)}(\dual K)$, that is $\d_0 p \in C^1(\dual K)$ is a
dual 1-cochain placed on the dual edges. This finally leads to the
conclusion that discrete pressure is a dual 0-cochain, that is, $p \in
C^0(\dual K)$ and thus the pressure must be placed at the
circumcenters of the top dimensional simplices.

Since $k$ is a dual 0-cochain we must use the discrete operator
$\d_0^\ast$ to replace the exterior derivative in
equation~\eqref{eq:dfdrcy}. Referring to the
diagrams~\eqref{eq:cochncmplx2d} and~\eqref{eq:cochncmplx3d} it is
clear that the discrete Hodge that should be used to replace $\hodge$
in equation~\eqref{eq:dfdrcy} is $(-1)^{n-1}\hodge^{-1}_{n-1}$, the
$(-1)^{n-1}$ sign coming from
expression~\eqref{eq:dscrtstrstr}. Finally, since $f\in C^{n-1}(K)$
clearly the discrete exterior derivative $\d_{n-1}$ will replace the
smooth $\d_{n-1}$ when discretizing equation~\eqref{eq:dfcnt}.

Thus the discretized equations corresponding to
equations~\eqref{eq:dfdrcy}-\eqref{eq:dfbc} are the very similar
looking
\begin{alignat}{2}
  f + \frac{\kappa}{\mu} \bigl((-1)^{n-1}\hodge_{n-1}^{-1} \,
  \d_0^\ast \,p\bigr) &= 0
 & &\quad \text{in } K \label{eq:darcy}\, ,\\
  \d_{n-1}\, f  &= \phi\omega 
  & &\quad \text{in } K \label{eq:cnt}\, ,\\
  f &= \psi\gamma 
  & &\quad\text{on } L \label{eq:bc}\, ,
\end{alignat}
where the unknowns and data are the cochains $f \in C^{n-1}(K)$, $p
\in C^0(\dual K)$, $\phi\omega \in C^n(K)$ and $\psi \gamma \in
C^{n-1}(K)$ where the last one is carried by $L$. The matrix
representation of equations~\eqref{eq:darcy} and~\eqref{eq:cnt}
adjusted for the boundary condition~\eqref{eq:bc} is the linear system
to be solved which is described next.

Let $f$ be the vector representing the cochain $f$ in the basis
\[
\bigl(\sigma_0^\ast, \ldots, \sigma_{N_{n-1}}^\ast\bigr)
\]
for $C^{n-1}(K)$ described in Section~\ref{subsec:dscrtfrms}. Recall
that $\sigma_i^\ast$ is the $(n-1)$-cochain that is 0 on $\sigma_i$,
the $(n-1)$ dimensional simplex number $i$. As mentioned in
Section~\ref{subsec:dscrtfrms} the simplices are ordered, for example,
in dictionary order \cite{BeHi2011}. Similarly, let $\kappa$,
$\phi\omega$ and $\psi\gamma$ be the vectors corresponding to the
other quantities appearing in
equations~\eqref{eq:darcy}-\eqref{eq:bc}.

To obtain the linear system to solve we first write
equations~\eqref{eq:darcy} and~\eqref{eq:cnt} in block matrix form
using the matrix form of the operators and objects. This yields
\begin{equation}
  \begin{bmatrix}
    I & (\kappa/\mu)(-1)^{(n-1)}\M_{n-1}^{-1} [\d_0^\ast]\\[1em] 
    \D_{n-1} & 0
  \end{bmatrix}
  \begin{bmatrix}
    f\\[1em]  p
  \end{bmatrix}  = 
  \begin{bmatrix}
    0\\[1em] \phi \omega
  \end{bmatrix} \, ,
  \label{eq:mtxwthI}
\end{equation}
where $I$ is an $N_{n-1}\times N_{n-1}$ identity matrix, 0 is an $N_n
\times N_n$ zero matrix and $[\d_0^\ast]$ is the matrix form of
$d_0$. Then using the fact that $[\d_0^\ast] = (-1)^nD^T_{n-1}$ we get
\begin{equation}
  \begin{bmatrix}
    I & -(\kappa/\mu)\M_{n-1}^{-1} \D_{n-1}^T\\[1em] 
    \D_{n-1} & 0
  \end{bmatrix}
  \begin{bmatrix}
    f\\[1em]  p
  \end{bmatrix}  = 
  \begin{bmatrix}
    0\\[1em] \phi \omega
  \end{bmatrix} \, .
  \label{eq:mtxwthIfnl}
\end{equation}

Assuming that the domain has only one connected component, the
pressure is unique only up to a constant. Hence, pressure at a single
point must be fixed to an arbitrary value to get a unique solution. To
impose the boundary conditions the right hand side of
equation~\eqref{eq:mtxwthIfnl} is adjusted for the known boundary
fluxes. Such an adjustment is also done for the assumed pressure at
one point. This is a standard procedure and it is described here
briefly for completeness. The adjustments are done by taking the
linear combination of the columns of the linear system matrix
in~\eqref{eq:mtxwthIfnl} corresponding to the known $f$ and $p$. The
coefficients in the linear combination are the known $f$ and $p$
values. The result is subtracted from the right hand side. These
columns and corresponding rows are then deleted from the matrices in
equation~\eqref{eq:mtxwthIfnl}.

Equation~\eqref{eq:mtxwthIfnl} can be written in a simpler standard
saddle point form \cite{BeGoLi2005} by starting from
equation~\eqref{eq:starfdrcy} instead of~\eqref{eq:dfdrcy}. This is
equivalent to multiplying the first block row of
~\eqref{eq:mtxwthIfnl} by $-(\mu/\kappa)\M_{n-1}$ and the
resulting system is
\begin{equation}
    \begin{bmatrix}
      -(\mu/k) \M_{n-1} & \D_{n-1}^T\\[1em]
      \D_{n-1} & 0
    \end{bmatrix}
    \begin{bmatrix}
      f\\[1em] p
    \end{bmatrix}  = 
    \begin{bmatrix}
      0\\[1em] \phi \omega
    \end{bmatrix}\, .
  \label{eq:mtxwthstr}
\end{equation}
To take into account the boundary conditions~\eqref{eq:bc} and the
assumed pressure at a point, the adjustments described above are
applied to equation~\eqref{eq:mtxwthstr}. The matrix on the left hand
side of~\eqref{eq:mtxwthstr} is a saddle type matrix of the form
\[
  \begin{bmatrix}
    A &B^T\\  B &0
  \end{bmatrix} \, .
\]
Even after the boundary fluxes and assumed pressure are taken into
account the form of the matrix stays the same, using sub-matrices of
$A$ and $B$. Here the matrix $A$ is a diagonal matrix with nonzero
diagonal entries since $A = (-\mu/\kappa)\M_{n-1}$ and $\M_{n-1}$ is the
diagonal discrete Hodge matrix. One difference between our method and
Raviart-Thomas finite elements is that in the latter case $A$ is not a
diagonal matrix.

The system obtained by adjusting equations~\eqref{eq:mtxwthIfnl}
or~\eqref{eq:mtxwthstr} for boundary conditions and assumed pressure
at a point can be solved for the remaining pressures and the unknown
fluxes by using any of the standard techniques for saddle type
matrices \cite{BeGoLi2005}. For the results in this paper we used the
Schur complement reduction method \cite{BeGoLi2005} or a direct
  solve for the full matrix system. In our case Schur complement is
particularly simple since the $A$ matrix can be explicitly inverted
trivially.

\subsection{Flux visualization}
\label{subsec:fluxviz}
Once the flux has been determined, we visualize it by using Whitney
interpolation as described in Section~\ref{subsec:intrp} to get a
smooth $(n-1)$-form inside each $n$-simplex. We then obtain the
corresponding velocity vector fields sampled at barycenters. The
sampling could be done at any location or locations in the interior of
the $n$-simplices, not just at the barycenter. Given a flux value, we
can determine the velocity as follows. The flux $f$ is related to the
velocity by $f = \hodge (\v^\flat)$ which implies that $\hodge f =
\hodge \hodge \bigl(\v^\flat\bigr) = (-1)^{n-1} v^\flat$.

Consider first $n=2$ and let the value of the Whitney interpolated
flux 1-form at a sampling point be $a\, dx + b\, dy$ where $a$ and
$b$ are some constants. Then
\[
\v^\flat = -(\hodge f) = -\hodge(a\, dx + b\, dy) = 
b\, dx - a\, dy \, ,
\]
which implies that for standard metric in $\mathbb{R}^2$, the velocity
$v$ at that point is the vector $(b,-a)$. For $n = 3$ if the value of
the Whitney interpolated flux 2-form at a sampling point inside a
tetrahedron is $f = a\, dy \wedge dz + b\, dz \wedge dx + c\,
dx \wedge dy$ then the associated velocity is given by
\[
	\v^\flat = \hodge f =  \hodge 
        (a\, dy \wedge dz + b\, dz \wedge dx + c \, dx \wedge dy)
	= a\, dx + b\, dy + c\, dz \, ,
\]
which implies that for standard metric in $\mathbb{R}^3$, the velocity
$\v$ is the vector $(a,b,c)$.

\section{Heterogeneous Permeability and Hodge Star}
\label{sec:heterogeneous}

In many physical problems the Hodge star (which is an operator
depending on the metric) appears as a material dependent operator. For
example when Maxwell's equations are written in terms of differential
forms, the electric permittivity and the magnetic permeability are
both Hodge stars \cite{Bossavit1998e}. The permittivity Hodge star
relates the electric field 1-form to the electric induction 2-form and
the permeability Hodge star relates the magnetic flux density 2-form
to the magnetic field 1-form.

In this section we rewrite the flux form of Darcy's law (equation
\eqref{eq:dfdrcy}) in a form that permits its discretization when the
permeability is a spatially dependent scalar quantity. The scalar
permeability value is discretized as constant in each $n$-simplex.
From these scalar values of the permeability, a spatially dependent
discrete Hodge star operator is constructed. This involves combining
the scalar permeabilities across an $(n-1)$-simplex in a weighted
average that is suggested by DEC and described in this section. The
resulting matrix for this heterogeneous Hodge star is still a diagonal
matrix. We start with the equations \eqref{eq:dfdrcy} --
\eqref{eq:dfbc} with the exception that the term $(\kappa/\mu)(\hodge
\d p)$ is replaced by $(1/\mu)(\hodge^\kappa \d p)$. Thus Darcy
law~\eqref{eq:dfdrcy} now becomes
\begin{equation}
  \label{eq:heterogeneous_darcy_law}
  f + \frac{1}{\mu}\bigl(\hodge^\kappa \d p\bigr) = 0 \quad\text{in}\; M
  \, ,
\end{equation}
where the Hodge star operator $\hodge^\kappa$ is the spatially
dependent Hodge star that depends on the permeability of the
medium. The continuity equation and boundary condition stay the same
as equations \eqref{eq:dfcnt} and \eqref{eq:dfbc}.  This heterogeneous
Hodge star is discretized as a diagonal matrix where the diagonal
entry corresponding to an $(n-1)$-simplex $\sigma^{n-1}$ in the matrix
$(\M^k_{n-1})^{-1}$ is
\begin{equation}
  \frac{\vol{\sigma^{n-1}}}{\vol{\dual \sigma^{n-1}}}\;
  \frac{\kappa_+\,\vol{\dual \sigma^{n-1} \cap \sigma_+^n} \, +\,
    \kappa_-\, \vol{\dual \sigma^{n-1} \cap \sigma_-^n}}
  {\vol{\dual \sigma^{n-1}}}\, .
\label{eq:hetero_hodge}
\end{equation}
Here $\sigma_+^n$ and $\sigma_-^n$ are the two simplices that contain
$\sigma^{n-1}$, and $\kappa_+$ and $\kappa_-$ are the permeabilities in these
$n$-simplices. The dual edge $\dual \sigma^{n-1}$ points from
$\sigma_-$ into $\sigma_+$. The expression $\vol{\dual
  \sigma^{n-1} \cap \sigma_+^n}$ stands for the length of the portion
of dual edge $\dual \sigma^{n-1}$ that lies in $\sigma_+^n$ etc. For
$\kappa_+=\kappa_- =\kappa$ expression~\eqref{eq:hetero_hodge} reduces to
\[
\kappa \; \frac{\vol{\sigma^{n-1}}}{\vol{\dual \sigma^{n-1}}} \, ,
\]
which is the corresponding diagonal entry of
$\M_{n-1}^{-1}$, thus yielding the usual discretization of
the homogeneous Hodge star scaled by $\kappa$.

The ratio on the right in expression~\eqref{eq:hetero_hodge} can be
interpreted as a weighted average of permeabilities. Note that it is
not a simple arithmetic or geometric mean of the permeabilities. The
weights are the same as the ones that were used for averaging
piecewise constant vector fields along a shared face
in~\cite{Hirani2003}. In~\cite{Hirani2003} it was shown that these are
the unique weights that yield a discrete divergence theorem. See
\cite[Figure~5.4, Section~5.5 and Section~6.1]{Hirani2003} for more
details.

\section{Numerical Results}
\label{sec:nmrclrslts}
We illustrate the performance of the proposed DEC based numerical
method for Darcy flow using many standard test problems. In all the
figures here that show a velocity vector field, the flux $f$ has been
visualized as the corresponding vector field. As described in
Section~\ref{subsec:fluxviz} the velocity vector field $v$ is obtained
from $f$ by using the relationship $v^\flat = (-1)^{n-1}\hodge f$. The
velocity vector field is sampled at barycenters of top dimensional
simplices and displayed as arrows based at barycenters. The pressure
in most figures is displayed by plotting it against the $x$ coordinate
of the circumcenter which is where the pressures are defined. The only
exceptions are the pressure plots in Figure~\ref{fig:sqr336coscos} in
which the pressure is displayed by coloring the triangle with a single
color based on pressure value at the circumcenter.
Figures~\ref{fig:sqr8hxgn6cnstvel},~\ref{fig:sqr336cnstvel},
~\ref{fig:sqr336coscos},~\ref{fig:sqr55cnstvel},~\ref{fig:onefreev_16},
~\ref{fig:cube244cnstvel},~\ref{fig:flx_cnvrgnc},
and~\ref{fig:hmsphrhlDECcnvrgnc} were constructed from data generated
by a Python implementation of our numerical method that used the PyDEC
software~\cite{BeHi2011}; and
Figures~\ref{fig:adjcntsqrs16cnstvel},~\ref{fig:lyrdrect1480vel}
and~\ref{fig:lyrdrect1480pr}, were generated by a MATLAB
implementation of our numerical method. We are well aware that
numerical tests of convergence of a method do not constitute a
convergence proof. Theoretical convergence analysis of this method is
a topic of future research.

\subsection{Patch tests in 2D}
\label{subsec:patch2d}
The first results shown in Figure~\ref{fig:sqr8hxgn6cnstvel} are for
patch tests \cite{IrLo1983}. It is desirable that for simple meshes a
numerical method for Darcy flow should reproduce constant velocity and
linear pressure exactly up to machine precision. The boundary
condition in these tests is derived from a constant horizontal
velocity $(1,0)$. Thus, for example, in the square domain $\v \cdot \n
= \psi = -1$ on the left edge, $\psi = 1$ on the right edge, and 0 on
the top and bottom edges of the square. Keeping in mind the
orientations and the definition of $\gamma$, when the discretized
equations~\eqref{eq:darcy}-\eqref{eq:bc} are used, $f = 1$ on the left
and right edges and 0 on the top and bottom edges of the square.

We now explain the sign of $f$ in more detail. Suppose the bottom left
and top left corner vertices of the square domain in
Figure~\ref{fig:sqr8hxgn6cnstvel} are labeled $v_0$ and $v_1$ and the
edge between them is oriented from $v_0$ to $v_1$. We will use the
name $\sigma$ for this oriented edge $[v_0, v_1]$ and denote the
vector from $v_0$ to $v_1$ by $\vec{\sigma}$. Assume also that the
square, and hence the triangle to which $\sigma$ belongs, is oriented
counterclockwise. We use the same name $f$ for the 1-cochain $f$ of
equations~\eqref{eq:darcy}-\eqref{eq:bc} and the 1-form $f$ of
equations~\eqref{eq:dfdrcy}-\eqref{eq:dfbc}, but to be more precise
the cochain $f$ should be referred to as $R(f)$ where $R$ is the
de~Rham map of Section~\ref{subsec:dscrtfrms}. The following
calculation explains why $\eval{R(f)}{\sigma} = +1$ in this setting.
\begin{align*}
\langle R(f), \sigma \rangle &= \int_\sigma \psi \gamma = 
\int_\sigma (\v \cdot \n) \gamma = (\v \cdot \n) \int_\sigma \gamma \\
&=(-1)\int_\sigma \gamma = -\gamma(\vec{\sigma}) = 
-\omega(\n,\vec{\sigma}) = +1 \, .
\end{align*}
Here the third equality follows from the fact that $\v \cdot \n$ is
constant along $\sigma$, the fifth equality is true because $\sigma$
is a straight line, and $\omega(\n,\vec{\sigma}) = -1$ because the
length of $\sigma$ is 1 and the basis $\bigl(\n,\vec{\sigma}\bigr)$
for the plane is oriented clockwise, which is opposite of the
orientation of the square. If $\sigma$ had been oriented from $v_1$ to
$v_0$ instead, the value of $\eval{R(f)}{\sigma}$ would have been $-1$
instead.

In this test, we are also given that $\phi = 0$, so that $\div v = 0$
(equivalently, $\d f = 0$) in the domain. The parameters $\kappa$ and $\mu$
are 1, and $g = 0$ so there is no external forcing. This example is
constructed by starting with pressure $p = -x + c$ for some constant
$c$. Then it follows that $v = (1,0)$ everywhere and $\div v = 0$ as
given. The numerical method is given the $\phi$ and $\psi$ and the
pressure and flux is computed using the method. Although such patch
tests do not guarantee that a method is high quality \cite{BaNa1997},
it is a convenient way to find problems with methods, as shown in
Figure~\ref{fig:sqr8nodal}.  If a method fails such a simple test, it
is probably unsuitable for the problem. 

The relative errors for the nonzero pressures are less than $3 \times
10^{-16}$ for the square and less than $7 \times 10^{-16}$ for the
hexagon shown in Figure~\ref{fig:sqr8hxgn6cnstvel}. Thus for these
simple meshes the relative error is close to machine precision. The
same test is repeated for a larger mesh of a square domain in
Figure~\ref{fig:sqr336cnstvel} for which the relative error in
pressure is less than $9 \times
10^{-12}$. Figure~\ref{fig:sqr55cnstvel} shows the result of same
patch test on a random Delaunay mesh, in which 37 vertices were placed
randomly, with no regard to mesh quality.

\subsection{Known solution with nonzero source term}
\label{subsec:general}
In the next test we compare the solution computed using our method
with an analytically known solution. Again the parameters $\kappa$ and
$\mu$ are 1 and $g=0$. The solution is constructed by starting with
pressure $p = \cos(\pi x) \cos(\pi y)$ from which an expression for
$v$ is derived. From $v$ one computes the divergence to derive that
the source/sink term is $\phi = 2 \pi^2 \cos(\pi x) \cos(\pi y)$. The
boundary data $\psi$ is constructed from the analytically computed
$v$. The numerical method is given $\phi$ and $\psi$ and used to
compute $v$ and $p$ in the domain from that.
Figure~\ref{fig:sqr336coscos} shows a comparison of the computed
pressure and velocity with the analytical solution.

\subsection{Discontinuous permeability}
\label{subsec:dscnt}

One of the new results of our approach is a Hodge star that allows the
discretization of the equations in the case of spatially varying
scalar permeability. The permeability is taken to be constant in each
$n$-simplex. The resulting heterogeneous diagonal Hodge star matrix
was defined in Section~\ref{sec:heterogeneous}. An important aspect of
any numerical method for Darcy flow is how well it can address
discontinuities in permeability. Our heterogeneous discrete Hodge star
allows us to test this aspect of our method.
Figure~\ref{fig:adjcntsqrs16cnstvel} shows the results of a patch test
with discontinuous permeability in two adjoining domains. The domain
is a rectangle in which the triangles in the left half are given a
permeability of $\kappa_1$ and those in the right half are given a
permeability of $\kappa_2$. The fluid flows in from the left and exits from
the right. The pressure should be a piecewise linear function whose
slope depends on the permeability and the velocity in the domain
should be constant. This is demonstrated in the results from our
method.

\subsection{Layered medium}
\label{subsec:layered}

Another common test in the Darcy flow literature is when the
discontinuities in the permeability vary across the flow rather than
along it. Such a medium is typically called a layered medium in which
the various layers are given different permeabilities. In our layered
medium computation we perform two tests. In one the permeabilities
alternate between 1 and 5 and in another they alternate between 1 and
10. The fluid comes into the domain from the left and exits from the
right as in our patch tests. The velocity should be horizontal and
constant in a layer. It should be larger in the low permeability
layers as computed by our method and shown in
Figure~\ref{fig:lyrdrect1480vel}. The correct pressure profile should
be a linear function of $x$ and this is seen in
Figure~\ref{fig:lyrdrect1480pr}.

\subsection{Patch tests in 3D}
\label{subsec:patch3d}

Many numerical formulations perform well in 2D, but their natural
extensions to 3D do not perform well. Herein we show that the proposed
formulation performs well even in 3D. To illustrate this we considered
two different computational domains, which are shown in
Figures~\ref{fig:onefreev_16} and~\ref{fig:cube244cnstvel}. The first
domain is a polyhedron with 16 tetrahedra and the second domain is a
cube with 244 tetrahedra. The analytical solution is constant velocity
along $x$ direction, and pressure linearly varying along $x$
direction. The obtained numerical results are plotted in
Figures~\ref{fig:onefreev_16} and~\ref{fig:cube244cnstvel} which shows
that the numerical method performed well. For example, the latter
figure shows that the relative error in pressure is less than $2\times
10^{-13}$.

\subsection{Circumcenter versus barycenter}

\begin{figure}[htb]
  \centering
  \includegraphics[scale=0.4,trim=0in 0in 0in 0in, clip]
  {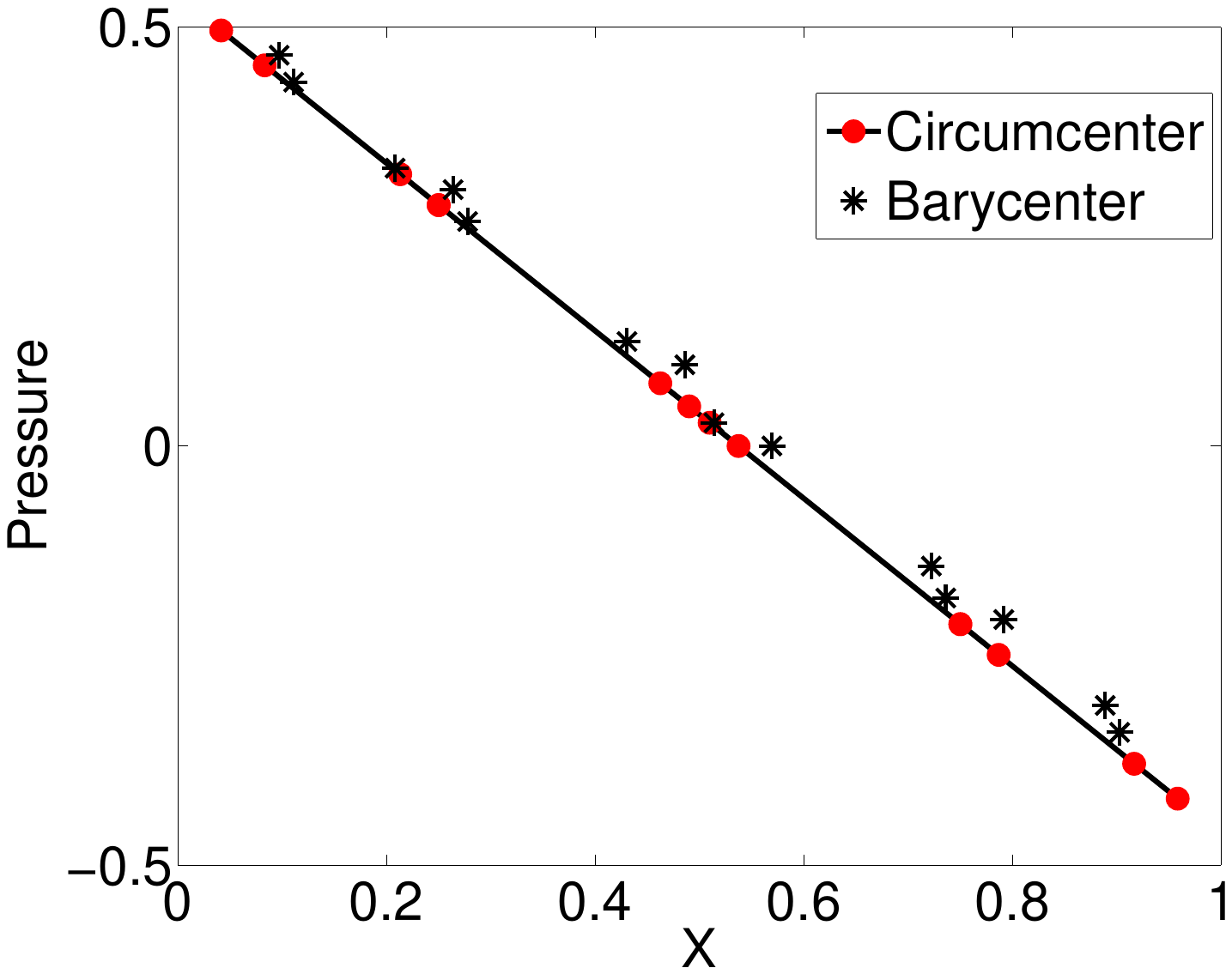}
  \caption{Shown here is the pressure computed for a square domain in
    which the boundary condition is given by constant horizontal
    velocity pointing to right. Thus velocity in the domain should be
    constant and pressure linear which is reproduced when pressures
    are located at circumcenters as dictated by DEC. See text for more
    discussion.}
 \label{fig:bary_vs_circum}
\end{figure}

\begin{figure}[p]
  \centering
  \subfigure{\includegraphics
    [scale=0.3, trim=0.5in 2.75in 0.5in 2.75in, clip]
    {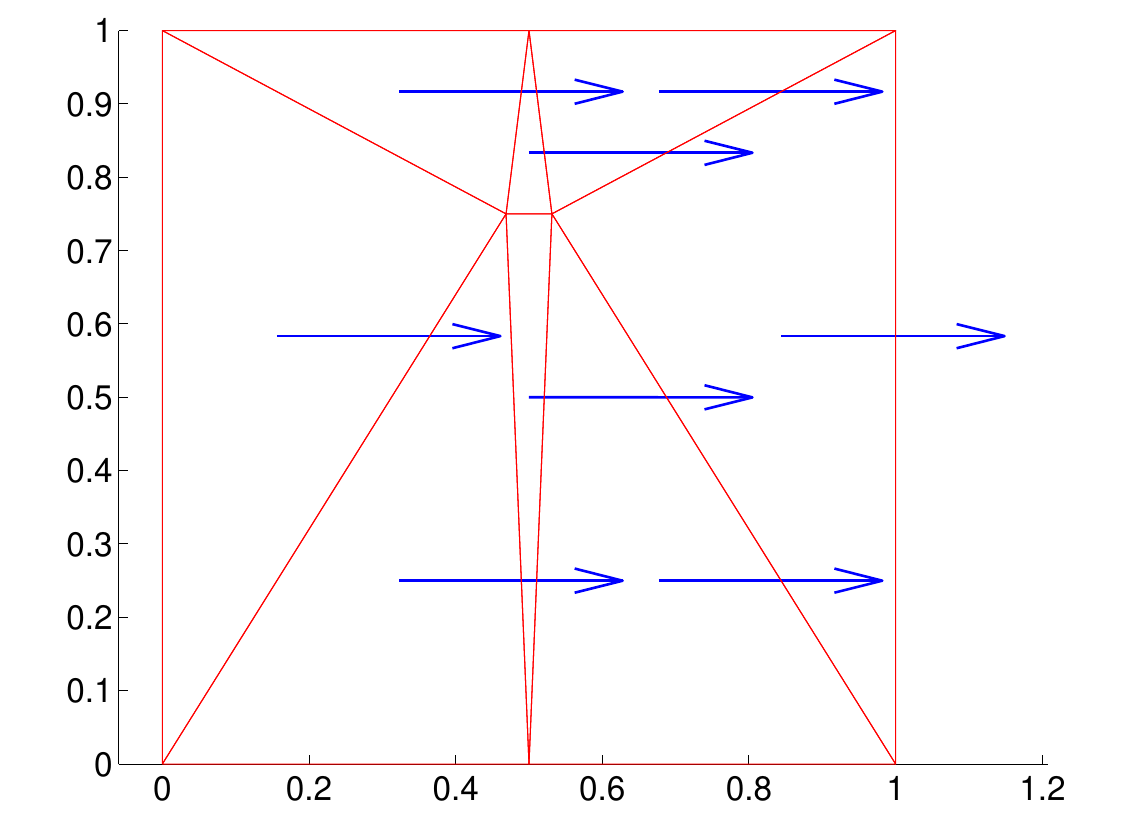}}
  \subfigure{\includegraphics
    [scale=0.3, trim=0.5in 2.5in 1in 2.75in, clip]
    {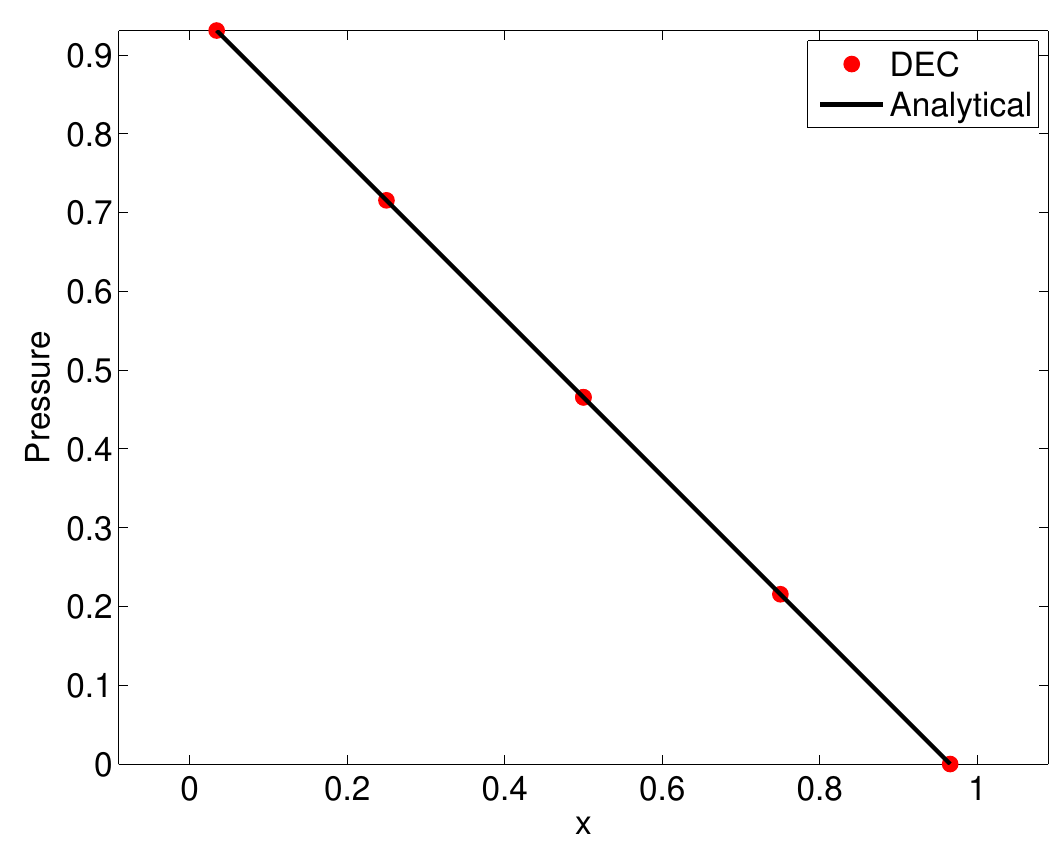}}
  \subfigure{\includegraphics
    [scale=0.3, trim=0.5in 2.75in 0.5in 2.75in, clip]
    {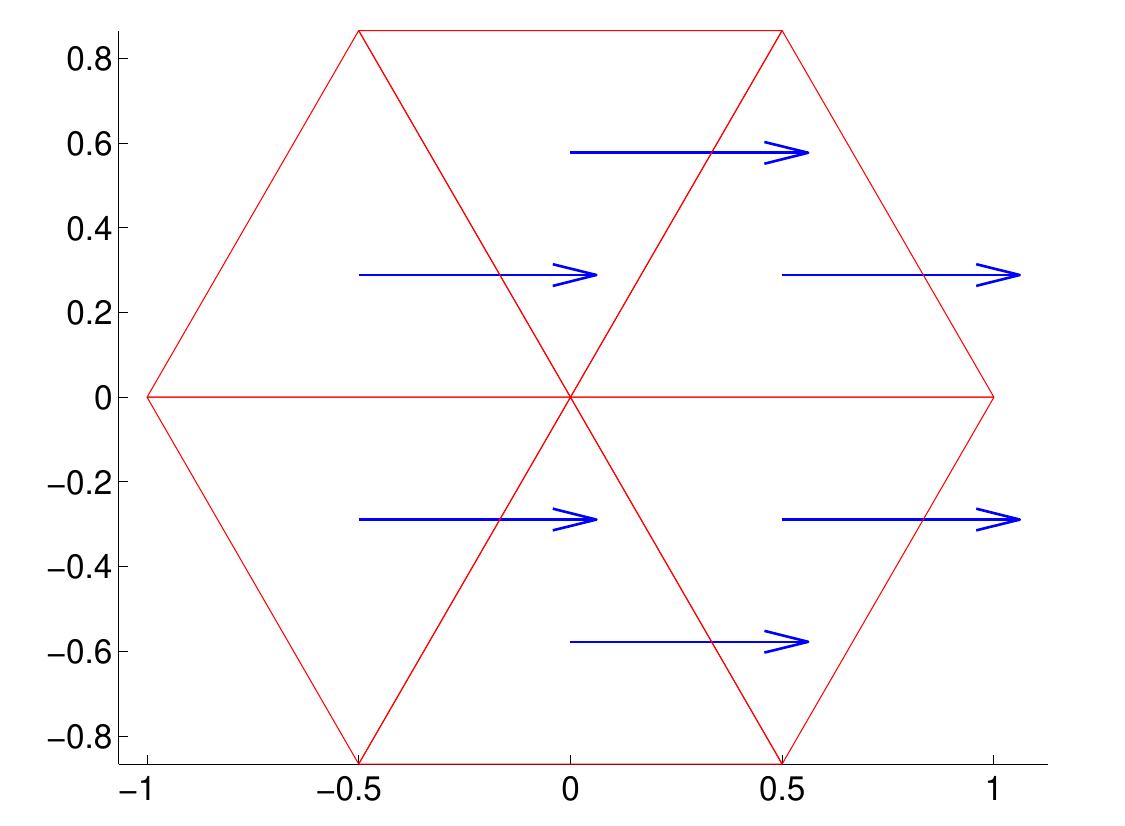}}
  \subfigure{\vspace{2in}}
  \subfigure{\includegraphics
    [scale=0.3, trim=0.5in 2.5in 1in 2.75in, clip]
    {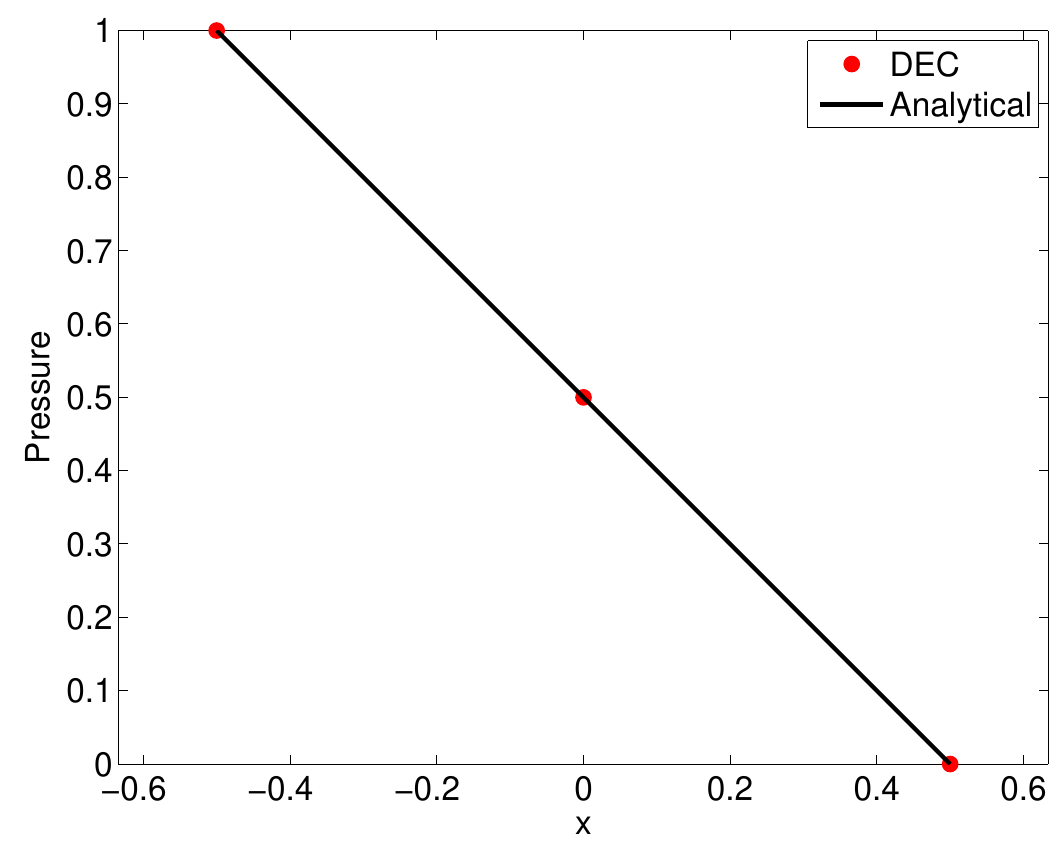}}
  \caption{Results from patch tests. The boundary conditions are
    obtained from a velocity of 1 in positive $x$ direction. Thus, for
    the square domain, in equation \eqref{eq:vfbc}, $\psi = -1$ on
    left edge, 1 on right ed 135900 ge, and 0 on top and bottom edges. For
    both the square and hexagon domains it is also given is that $\phi
    = 0$ in ~\eqref{eq:vfcnt}, constants $\kappa = 1$, $\mu = 1$ and
    external body force acceleration $\g = 0$ in~\eqref{eq:vfdrcy}.
    The pressure should be linear as shown in the right figures. The
    left figures show that the velocity is constant, as expected. The
    velocity displayed is interpolated from the flux through each
    edge, using Whitney 1-form and sampled and converted to a vector.}
  \label{fig:sqr8hxgn6cnstvel}
\end{figure}

\begin{figure}[ht]
  \centering
  \includegraphics[scale=0.5,trim=0.5in 2in 1in 2in,clip]
  {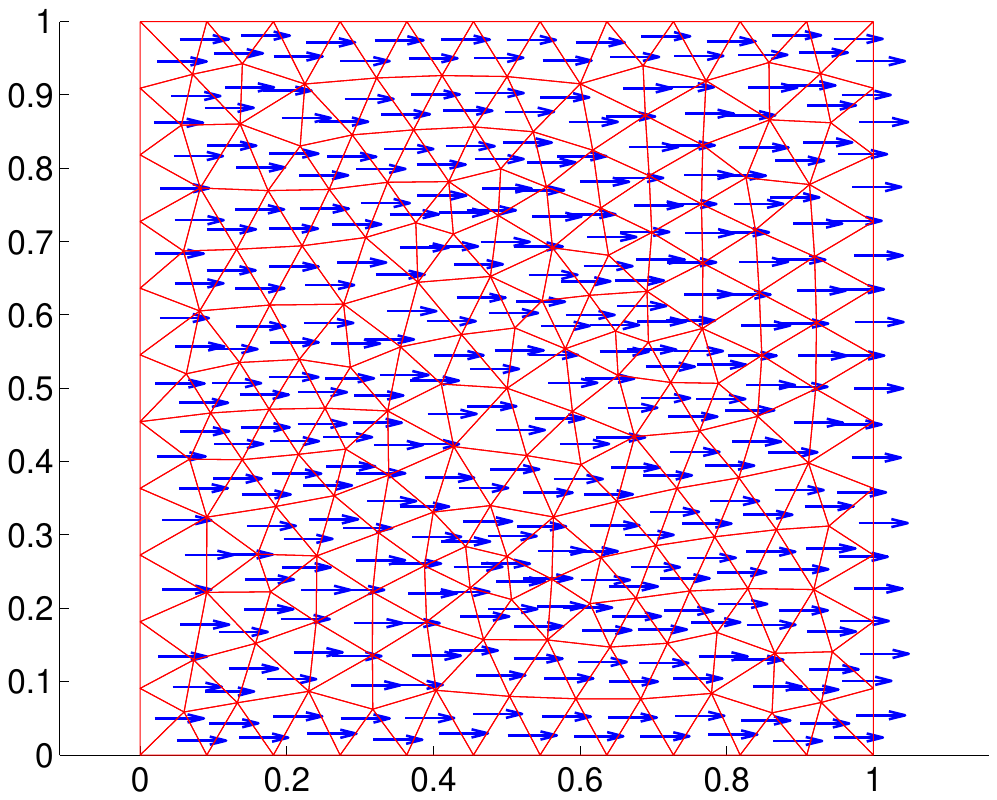}
 \caption{Patch test with larger mesh and with the same boundary
    conditions, parameters and other data as the square in
    Figure~\ref{fig:sqr8hxgn6cnstvel}.}
  \label{fig:sqr336cnstvel}
\end{figure}

\begin{figure}[hb]
  \centering
  \includegraphics[scale=0.5]
  {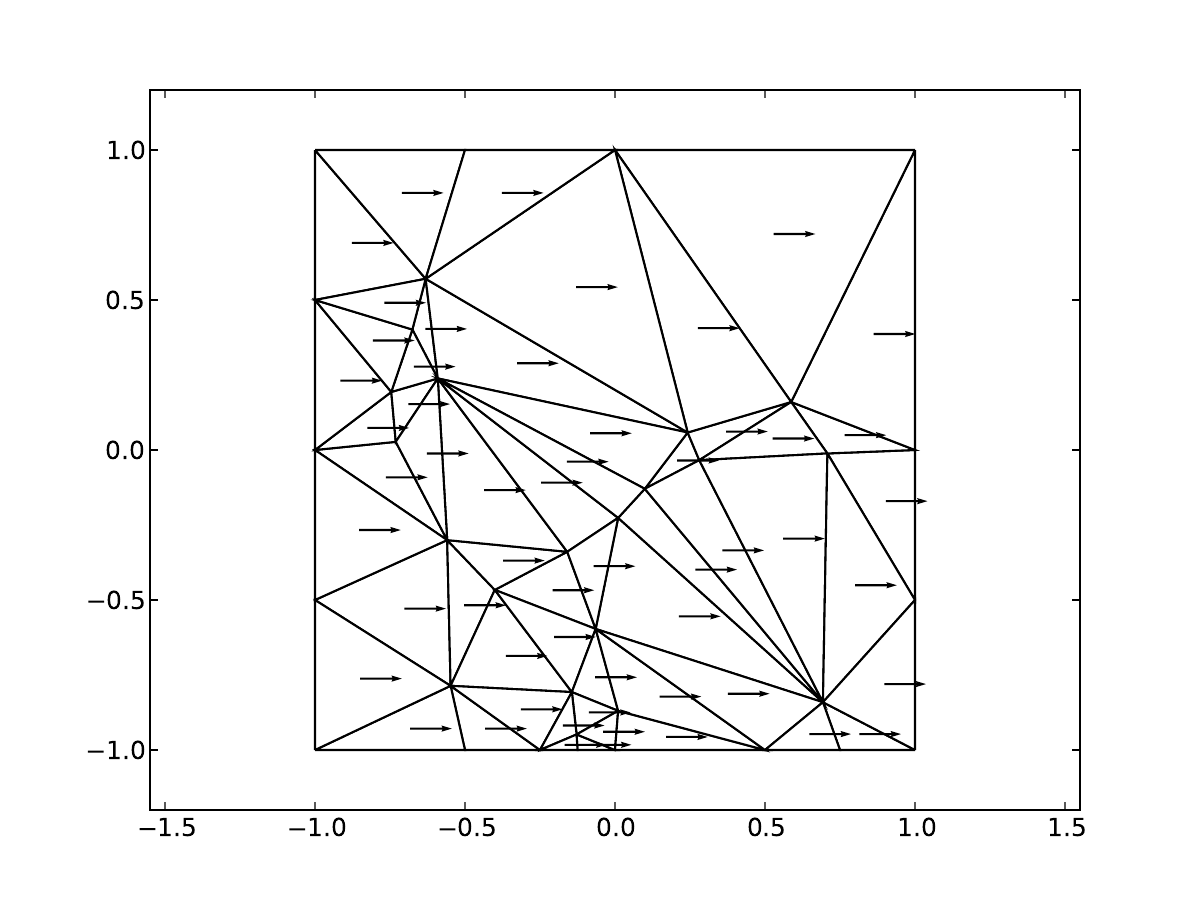}
  \includegraphics[scale=0.5]
  {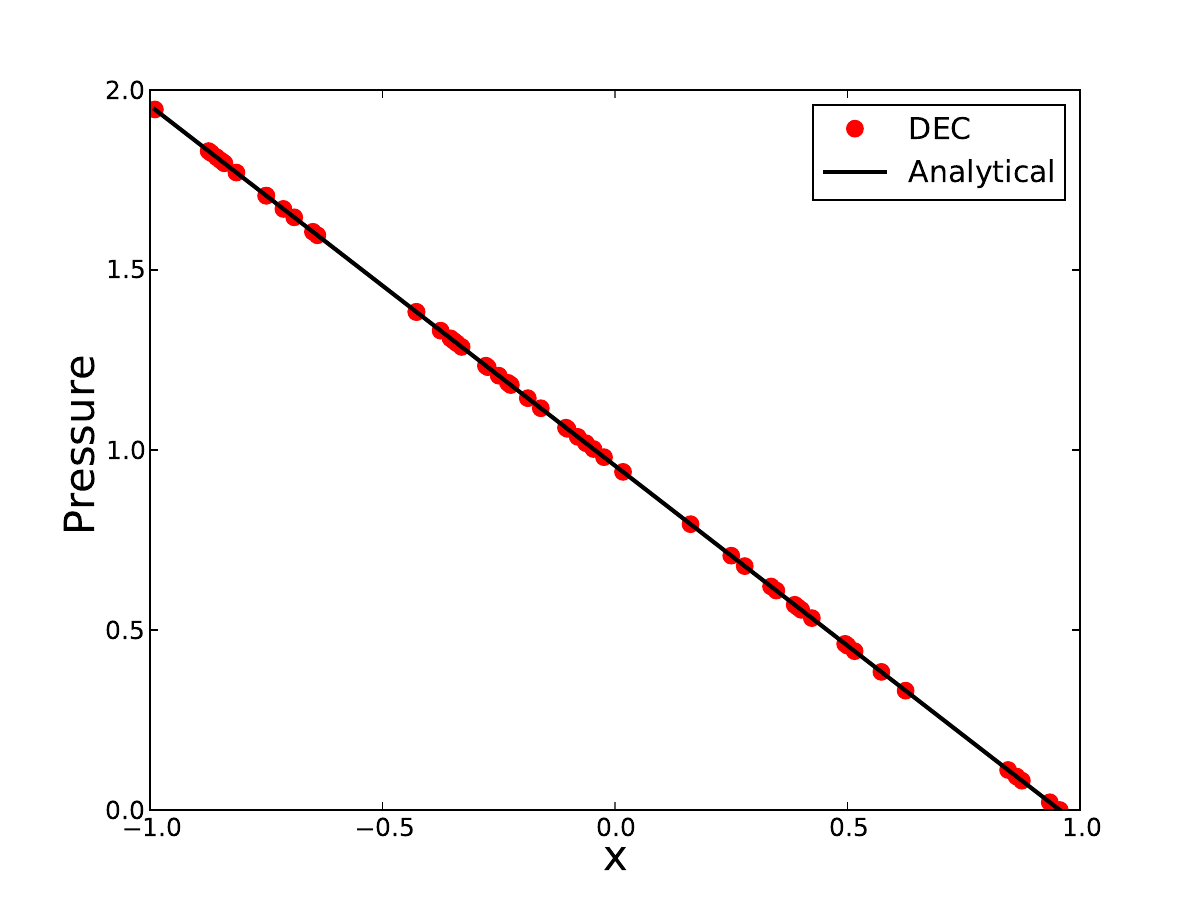}
  \caption{Patch test using a Delaunay triangulation of randomly
    placed 37 vertices in a square. Several triangles are obtuse
    angled and no attempt has been made to make the mesh have high
    quality. The boundary conditions, parameters and other data are
    same as for the square in Figure~\ref{fig:sqr8hxgn6cnstvel}. Note
    that circumcenters of obtuse angled triangles fall outside those
    triangles but that this does not reduce the effectiveness of the
    method.}
  \label{fig:sqr55cnstvel}
\end{figure}

\begin{figure}[ht]
  \centering
  \subfigure{
    \centering
    \includegraphics[scale=0.3,trim=0.5in 2.5in 0.5in 2.75in,clip]
    {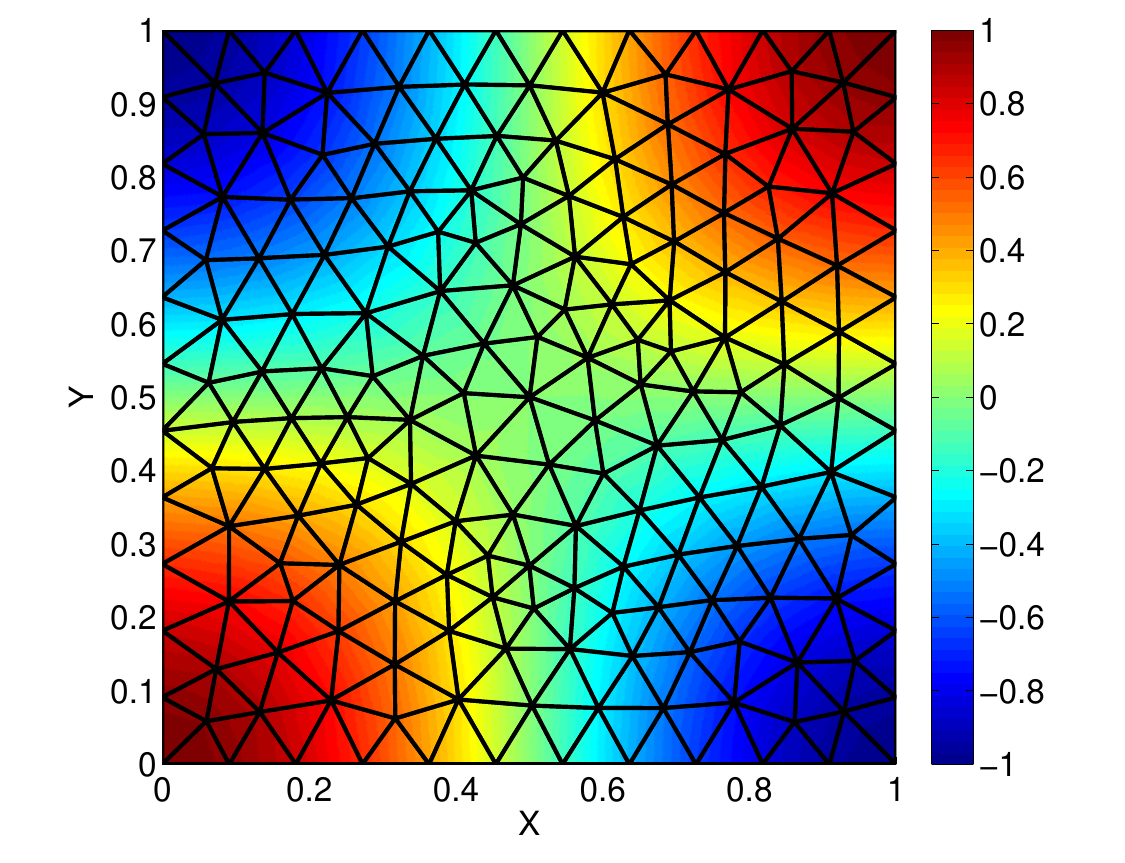}}
  \subfigure{
    \centering
    \includegraphics[scale=0.3,trim=0.5in 2.5in 0.5in 2.75in,clip]
    {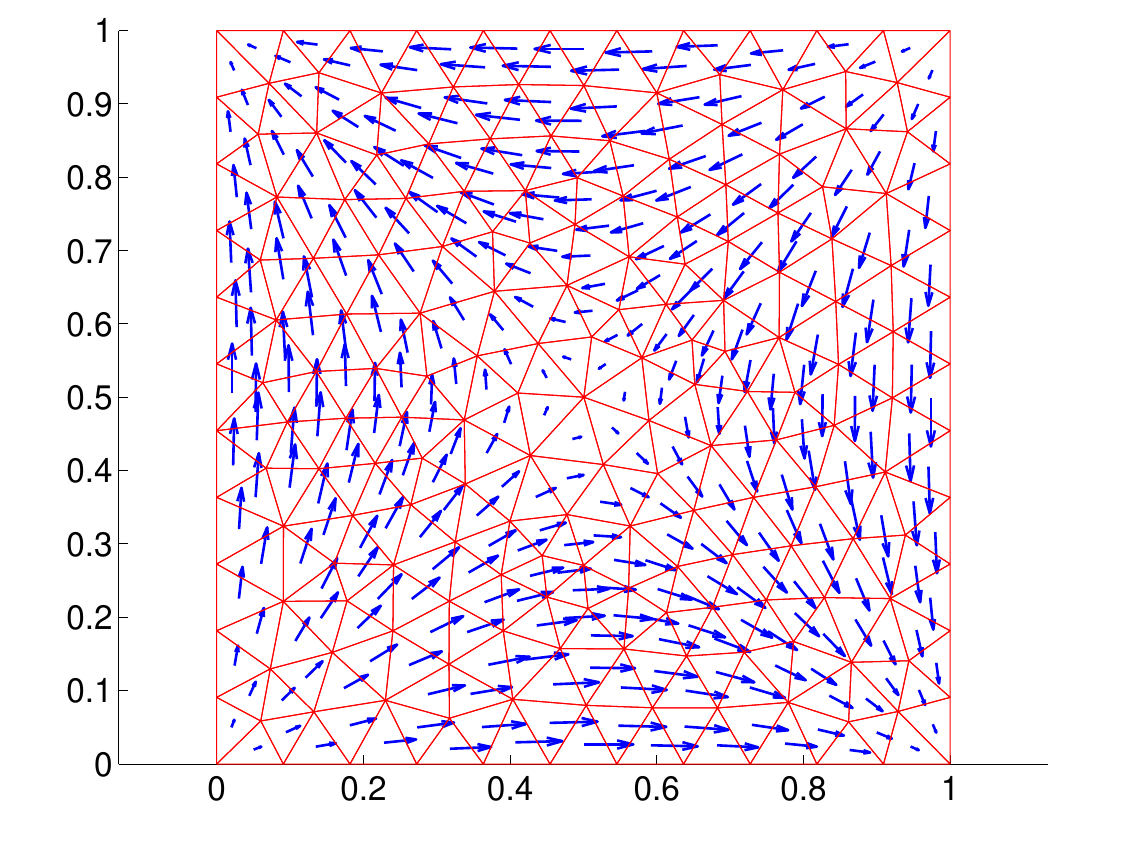}}
  \subfigure{
    \centering
    \includegraphics[scale=0.3,trim=0.5in 2.5in 0.5in 2.75in,clip]
    {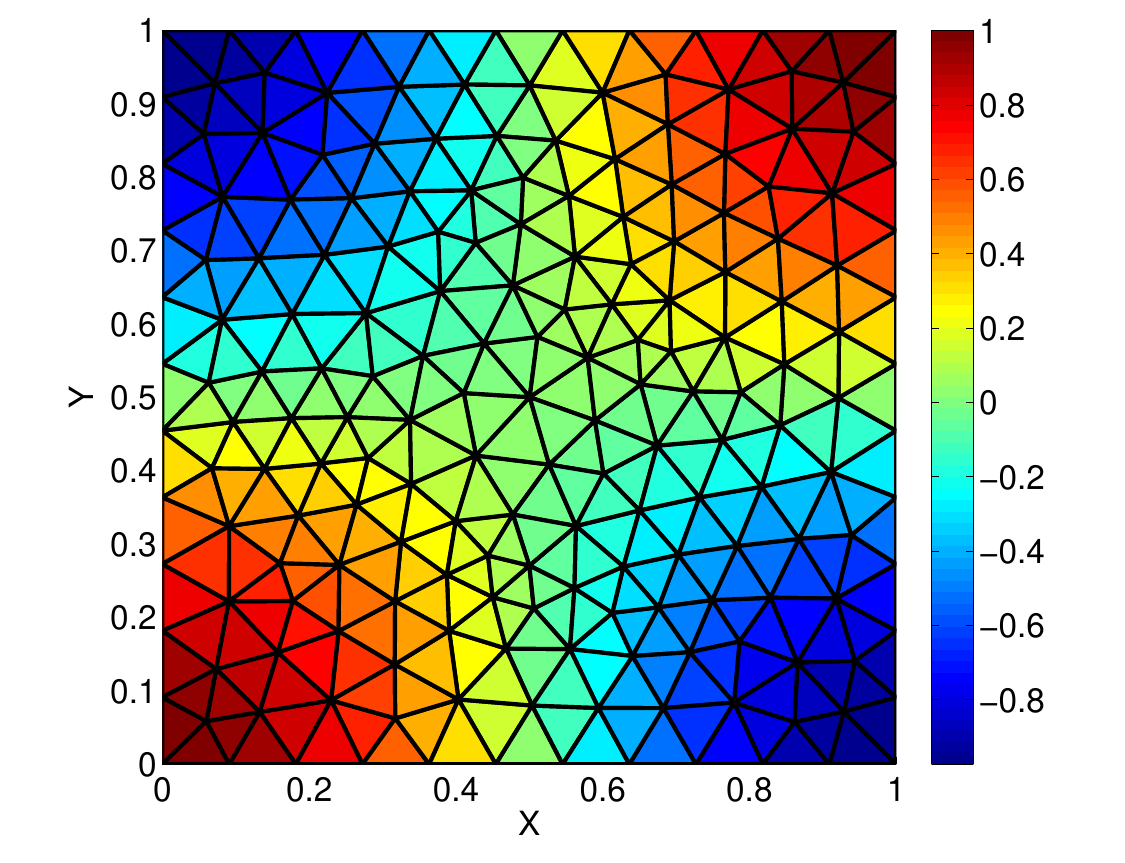}}
  \subfigure{
    \centering
    \includegraphics[scale=0.3,trim=0.5in 2.5in 0.5in 2.75in,clip]
    {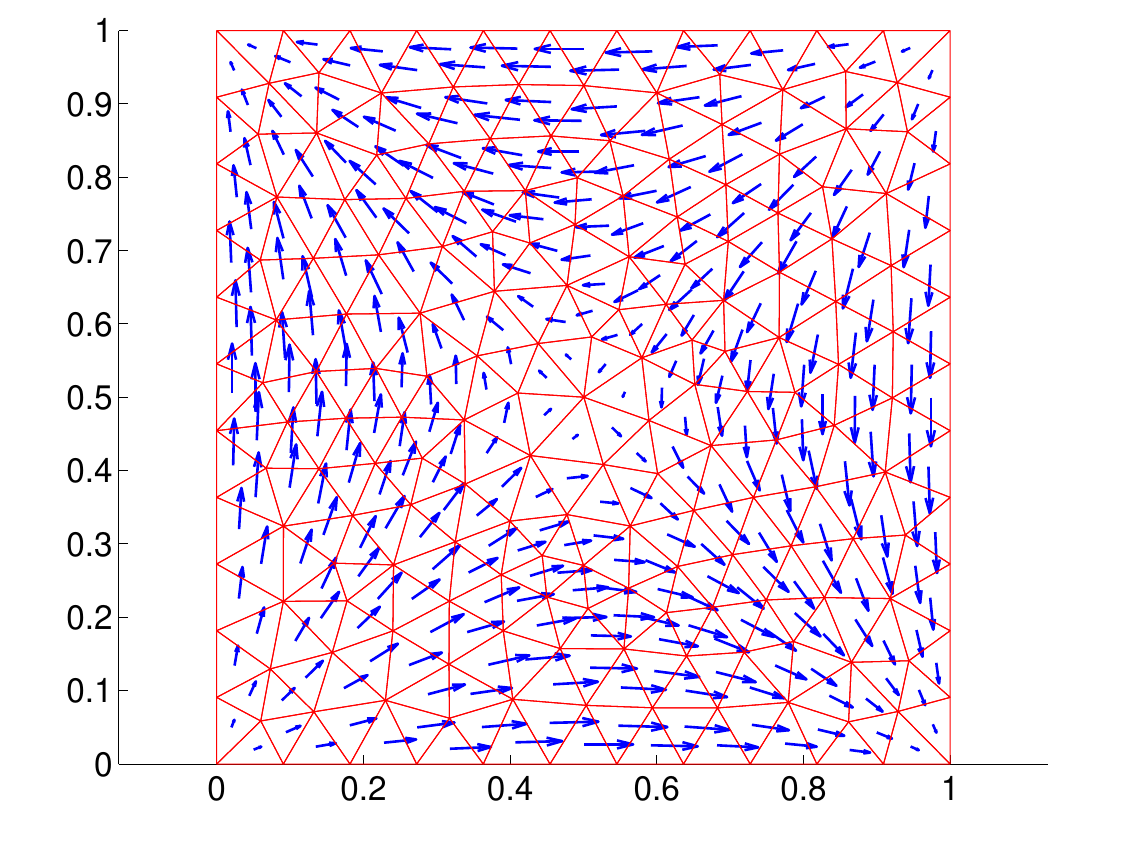}}

  \caption{Comparison of analytically known solution with a solution
    computed using the method developed here. The top left figure
    shows the pressure part of the analytical solution to the Darcy
    law problem with $\kappa$ = 1, $\mu$ = 1, $\g$ = 0, and $\phi$ = $2
    \pi^2 \cos(\pi x) \cos(\pi y)$. This implies $p = \cos(\pi x)
    \cos(\pi y)$ which is shown overlaid on the mesh. In the left
    bottom figure we show the pressure computed with the numerical
    method proposed here. Each triangle is colored by the pressure at
    its circumcenter. The top right figure shows the analytically
    computed vector field, and the numerically computed vector field
    is shown in the right bottom figure, which is obtained from the
    flux by interpolation of 1-forms.}
  \label{fig:sqr336coscos}
\end{figure}

\begin{figure}[ht]
  \centering
  \subfigure{\includegraphics[width=5.5in]
    {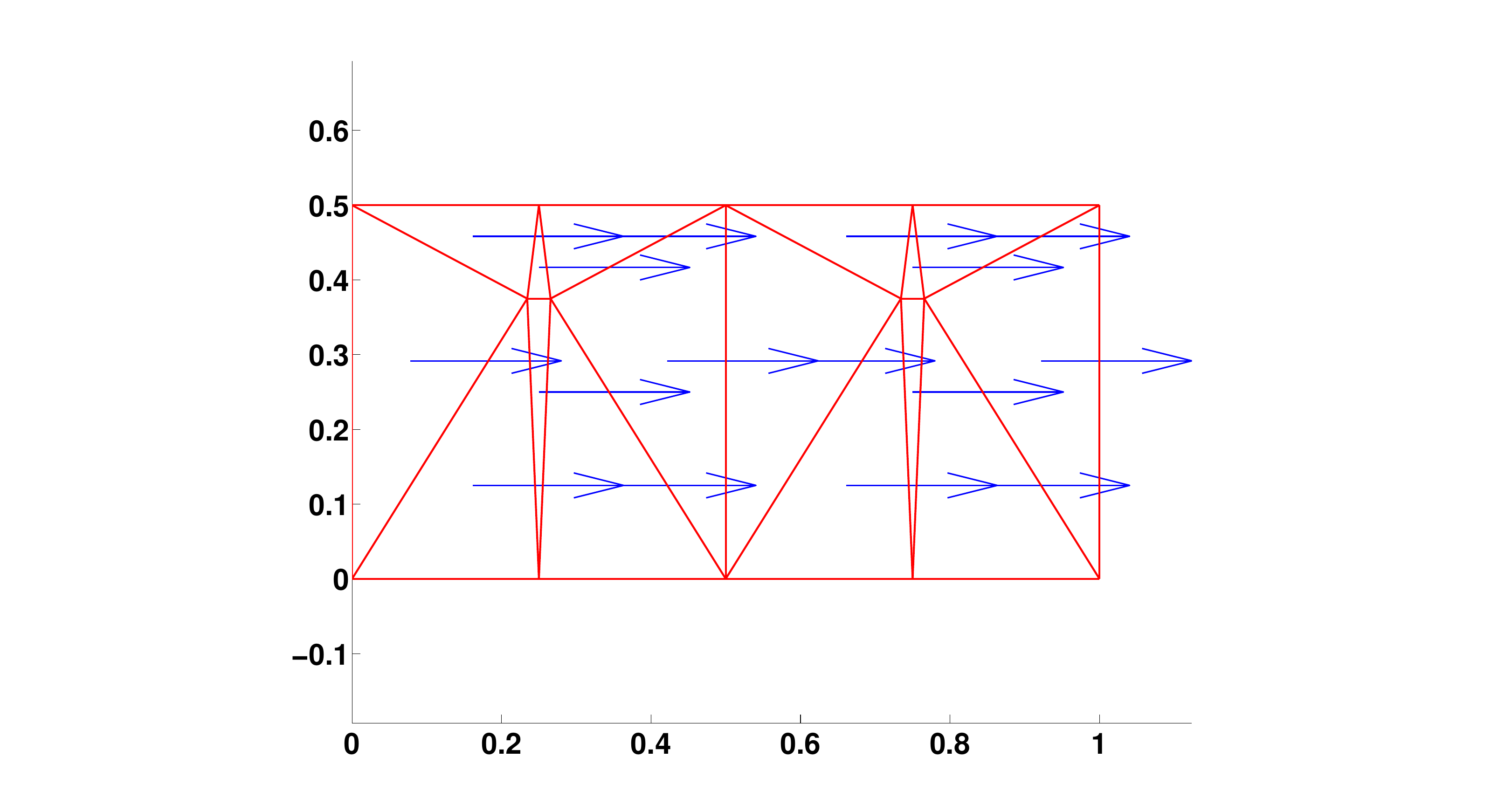}}
  \subfigure{\includegraphics[width=5.5in]
    {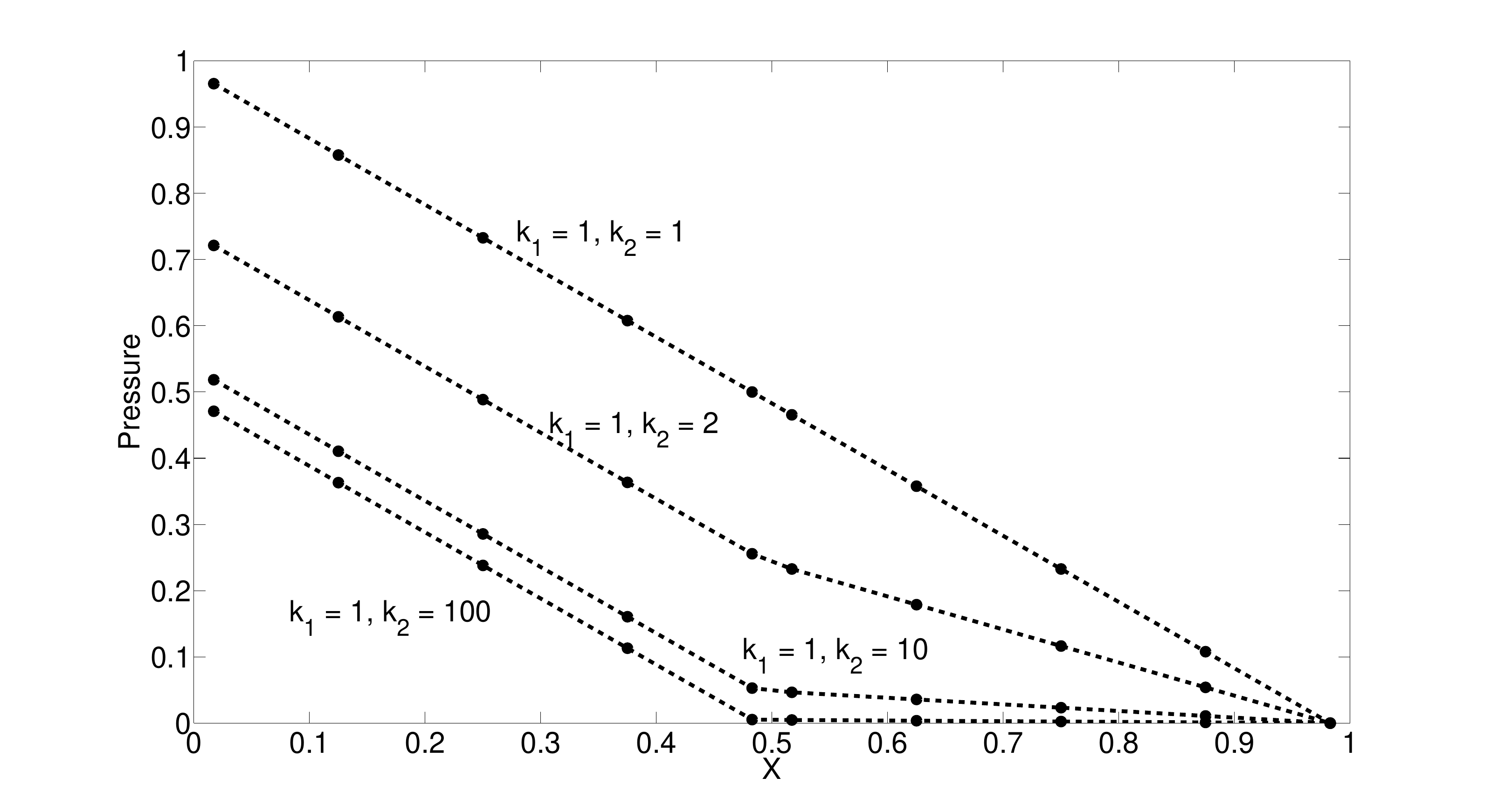}}
  \caption{In this test the left and right halves of the mesh are
    given different permeabilities, $\kappa_1$ on the left and
    $\kappa_2$ on right. Thus the permeability jumps across the middle
    vertical edge at $x = 0.5$. The values that we use for
    $(\kappa_1,\kappa_2)$ are $(1,1)$, $(1,2)$, $(1,10)$ and
    $(1,100)$.  The boundary condition is derived from velocity
    $(1,0)$ as in Figure~\ref{fig:sqr8hxgn6cnstvel}. The top figure
    shows the computed velocity interpolated from flux for
    $\kappa_1=1$ and $\kappa_2=10$. All other value pairs also result
    in a constant velocity solution as expected.  The pressure is
    piecewise linear with a jump at the discontinuity at $x=0.5$ as
    shown in the bottom figure.}
  \label{fig:adjcntsqrs16cnstvel}
\end{figure}

\begin{figure}[ht]
  \centering
  \subfigure{
    \includegraphics[scale=0.3, trim=2.5in 0.5in 0in 0.75in, clip]
    {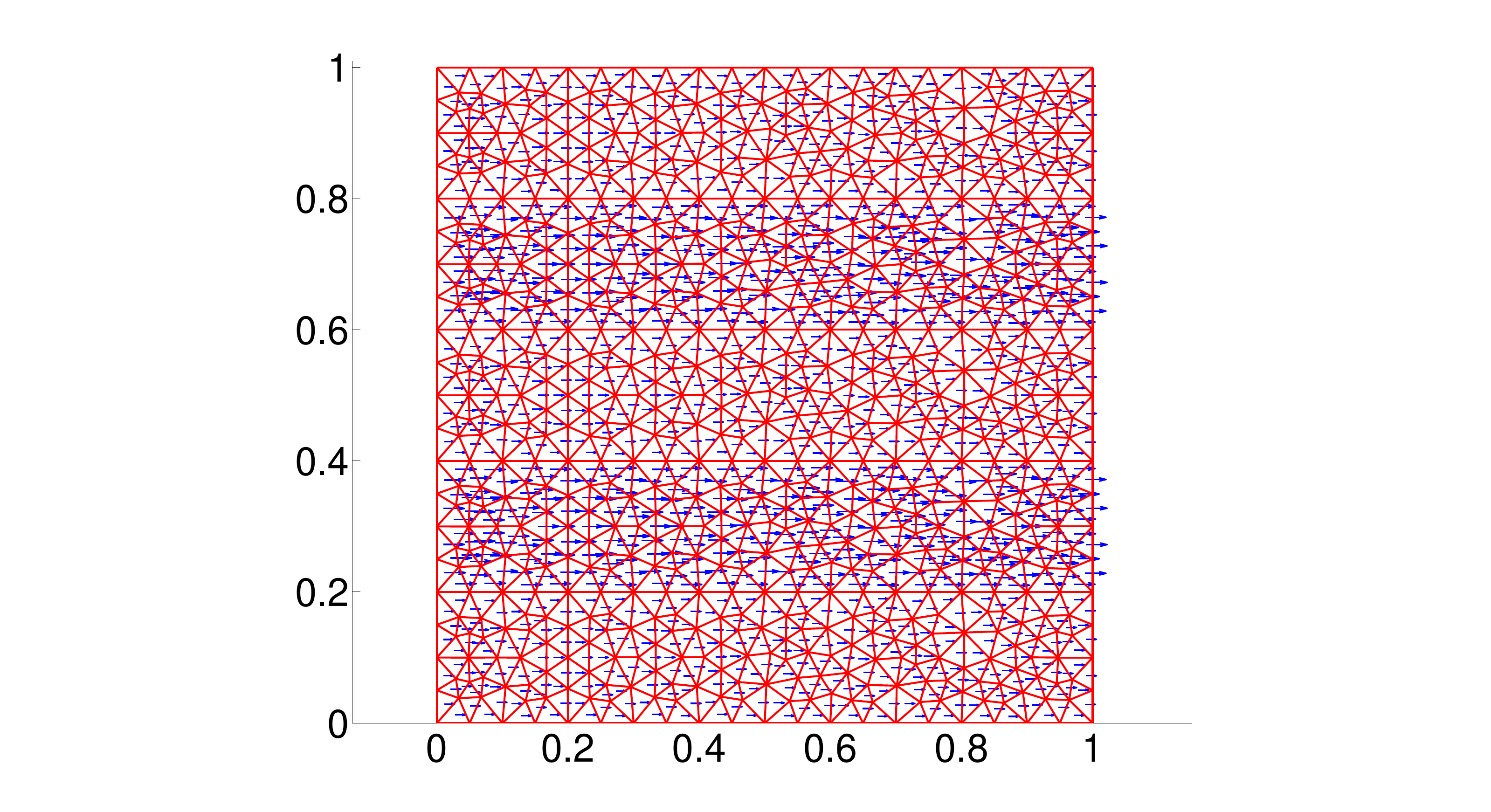}}
  \subfigure{
    \includegraphics[scale=0.3, trim=2.5in 0.5in 0in 0.75in, clip]
    {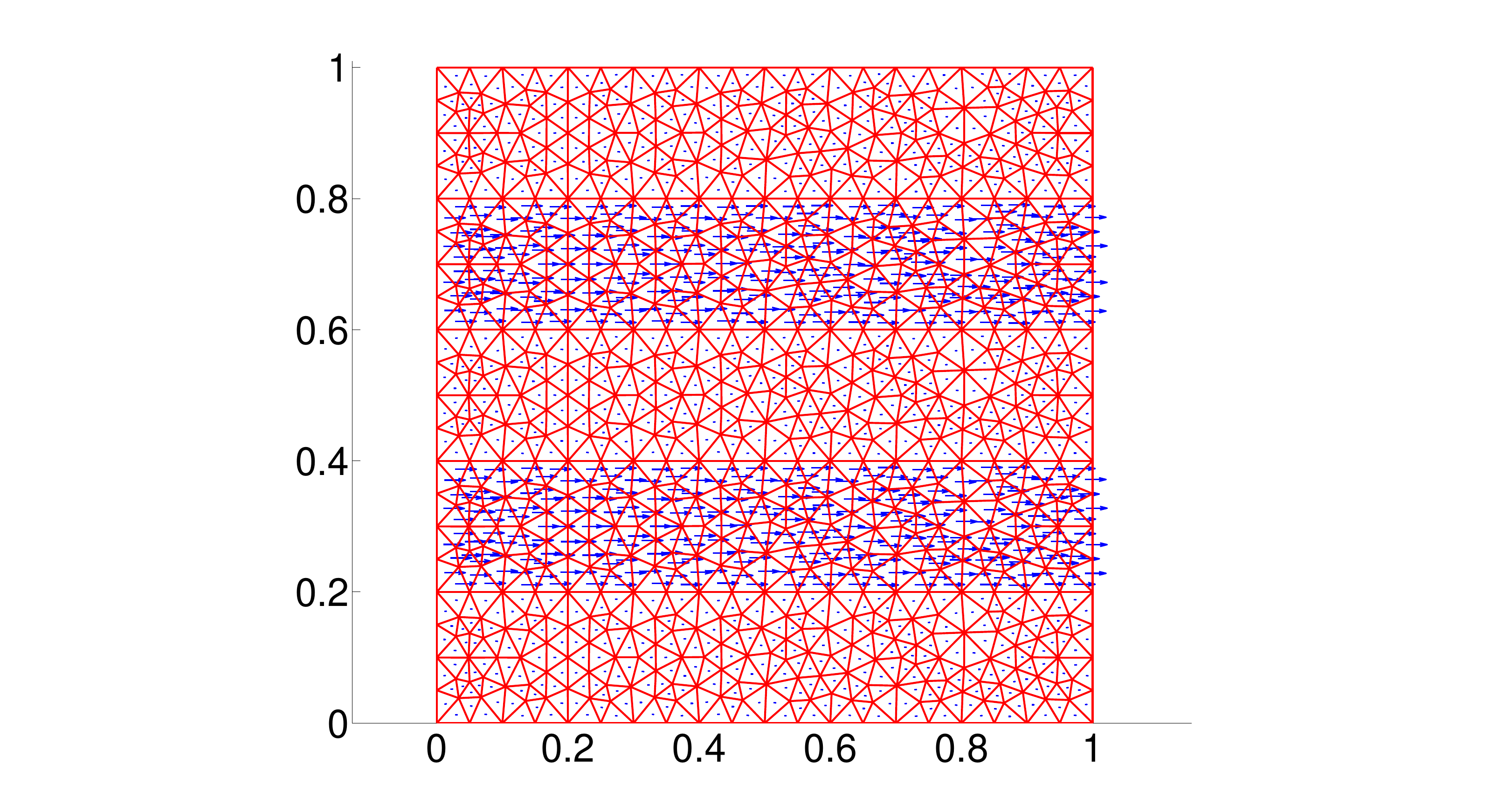}}
  \caption{Layered medium with 2 different permeability patterns. The
    domain has 5 layers with alternating permeability. In the top
    figure the permeability $\kappa$, from bottom layer to top is 5,
    10, 5, 10 and 5. In the bottom figure the permeability $\kappa$ is
    1, 10, 1, 10 and 1. The computed flux is visualized as a vector
    field.}
  \label{fig:lyrdrect1480vel}
\end{figure}

\begin{figure}[hb]
  \centering
  \subfigure{
    \includegraphics[width=5in]
    {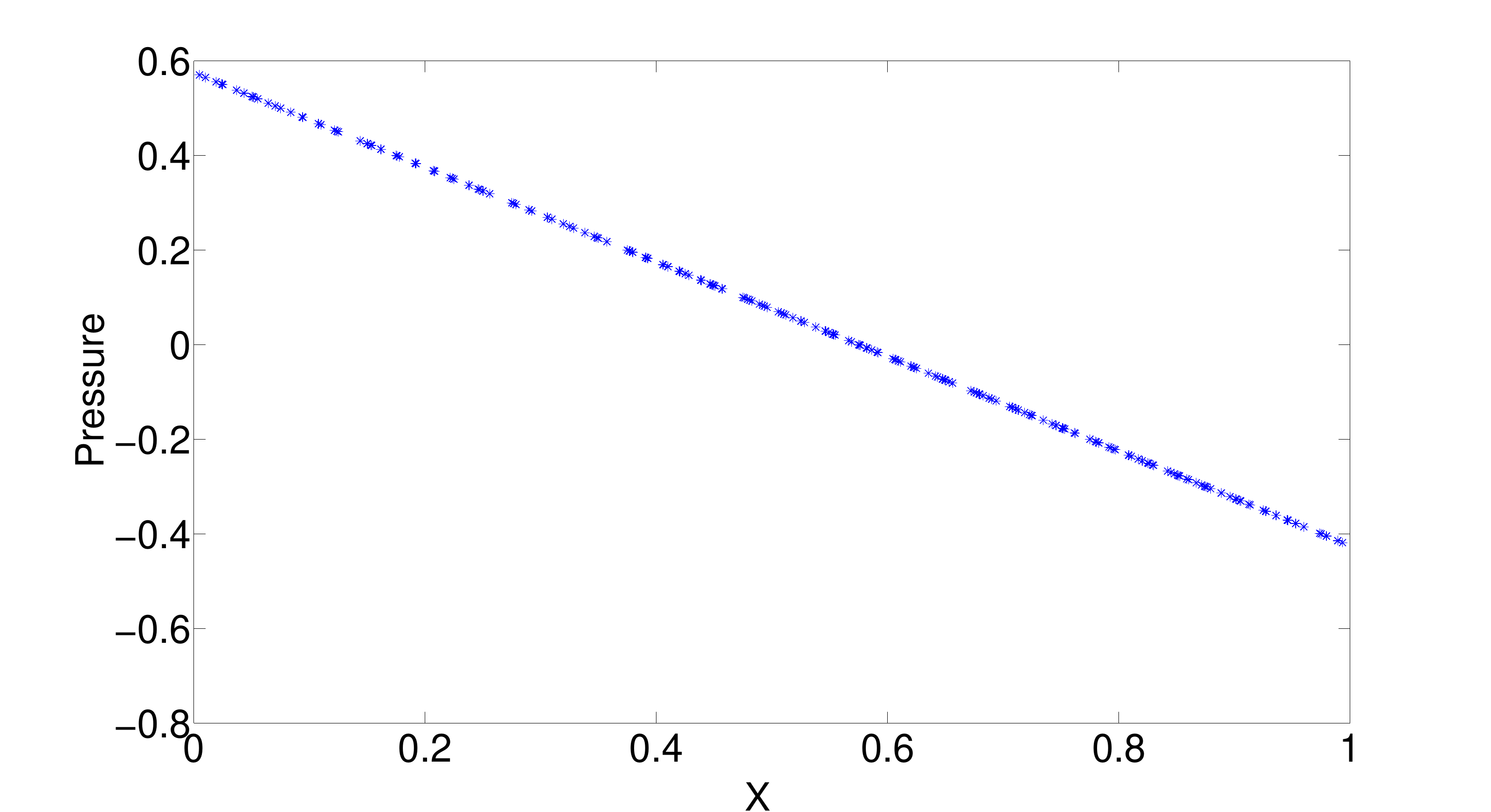}}
  \subfigure{
    \includegraphics[width=5in]
    {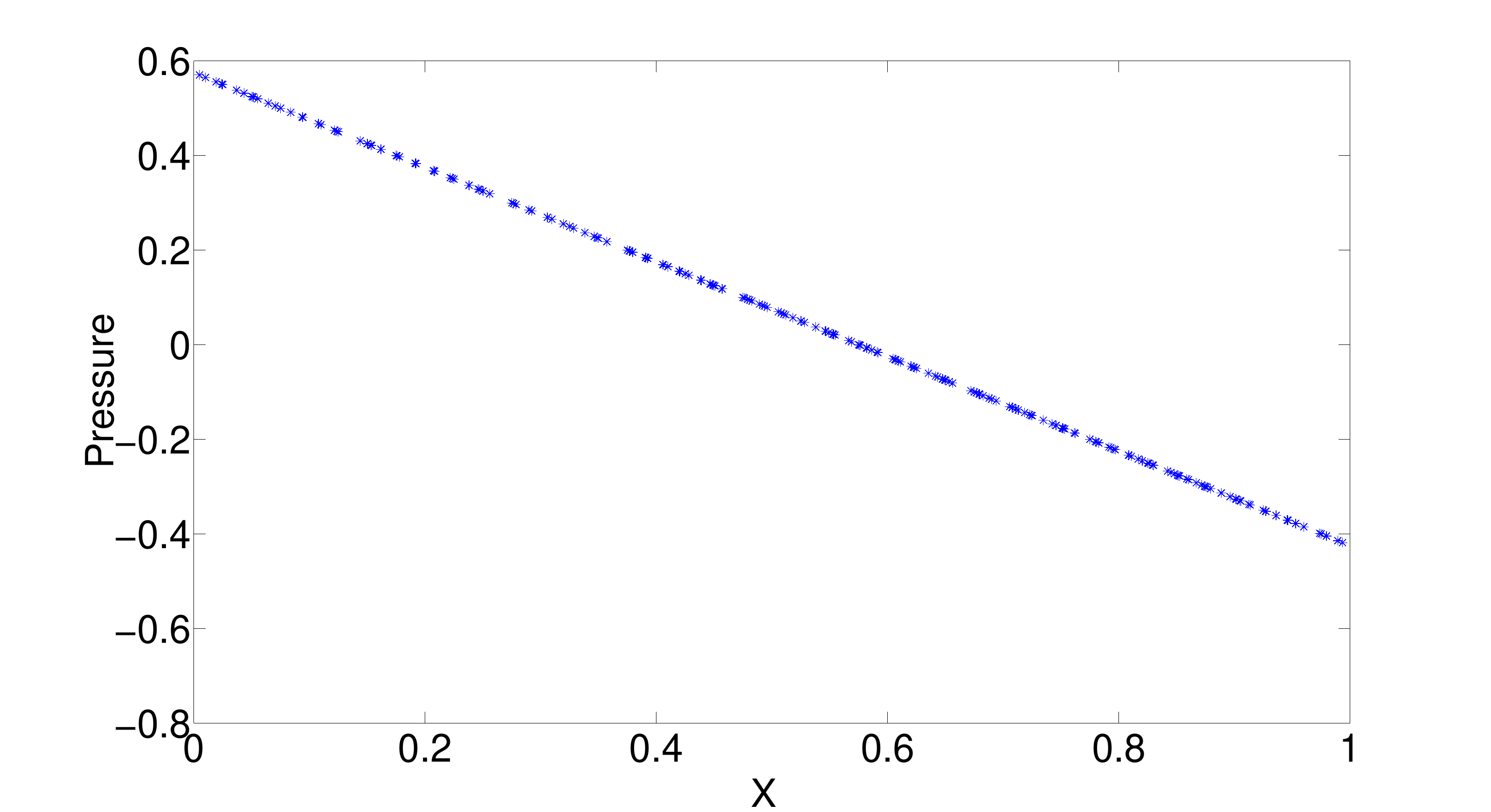}}
  \caption{Pressure for the layered medium show in
    Figure~\ref{fig:lyrdrect1480vel}. The pressure is linear as
    expected. The top and bottom figures correspond to the top and
    bottom figures in Figure~\ref{fig:lyrdrect1480vel}.}
  \label{fig:lyrdrect1480pr}
\end{figure}

\begin{figure}[ht]
  \centering
  \includegraphics[scale=0.4]{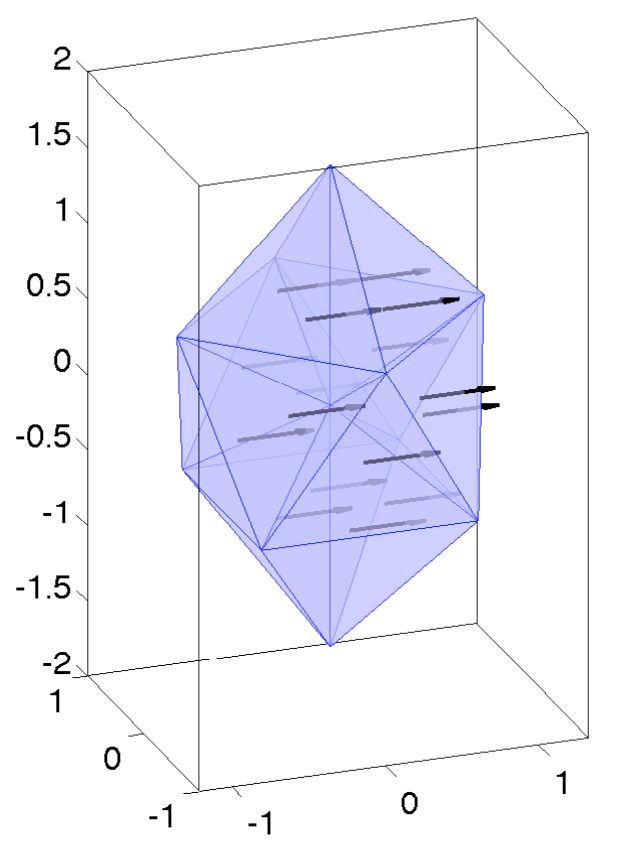}
  \includegraphics[scale=0.3]{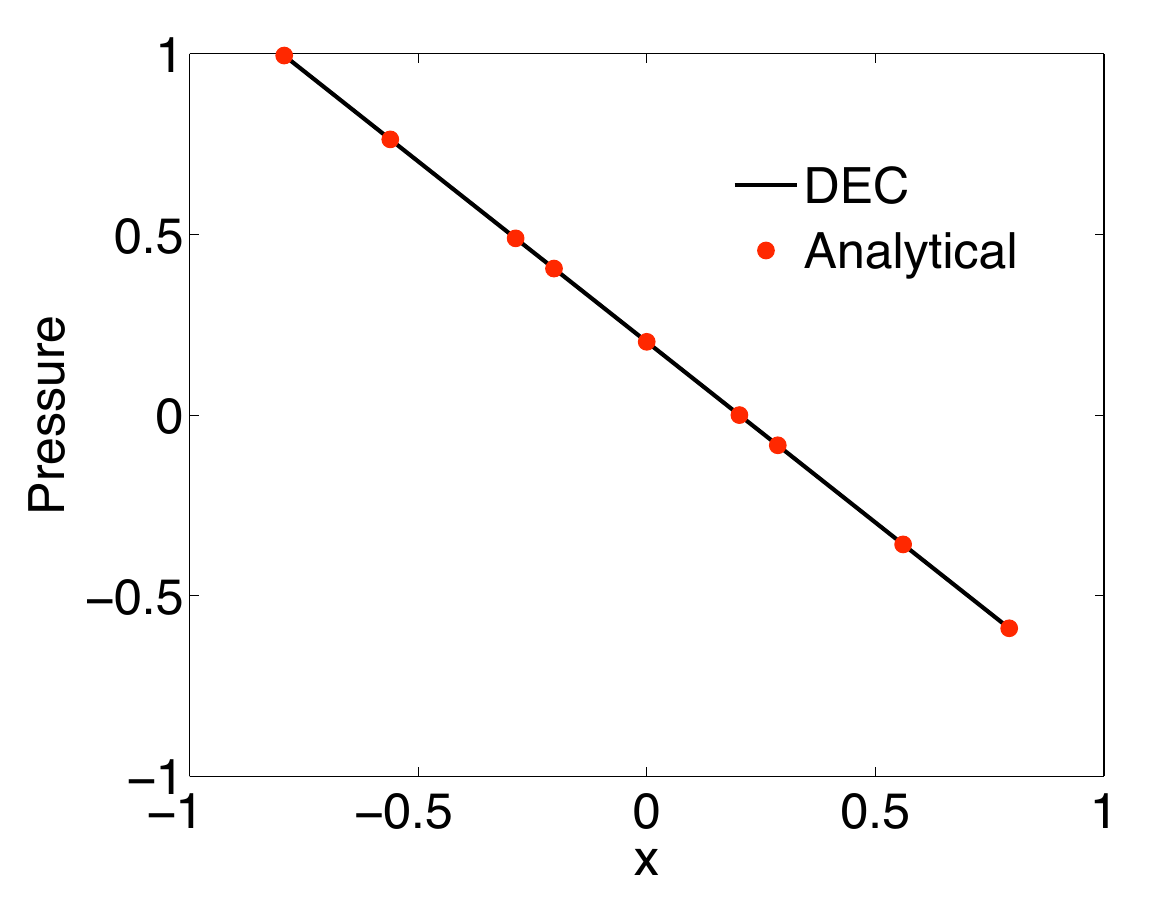}
  \caption{Patch test in 3D. The mesh shown has 16 tetrahedra.
    The fluid velocity on the boundary is from negative
    $x$ to positive $x$ direction, with no $y$ or $z$ components. The
    fluxes (2-cochains) across the internal faces are computed using
    our method and interpolated using the Whitney map.  This is then
    sampled at the circumcenters, converted into a vector field and
    plotted as arrows. }
  \label{fig:onefreev_16}
\end{figure}

\begin{figure}[ht]
  \centering
  \subfigure{
  \includegraphics[scale=0.4,trim=1in 2.5in 1.5in 2.5in,clip]
  {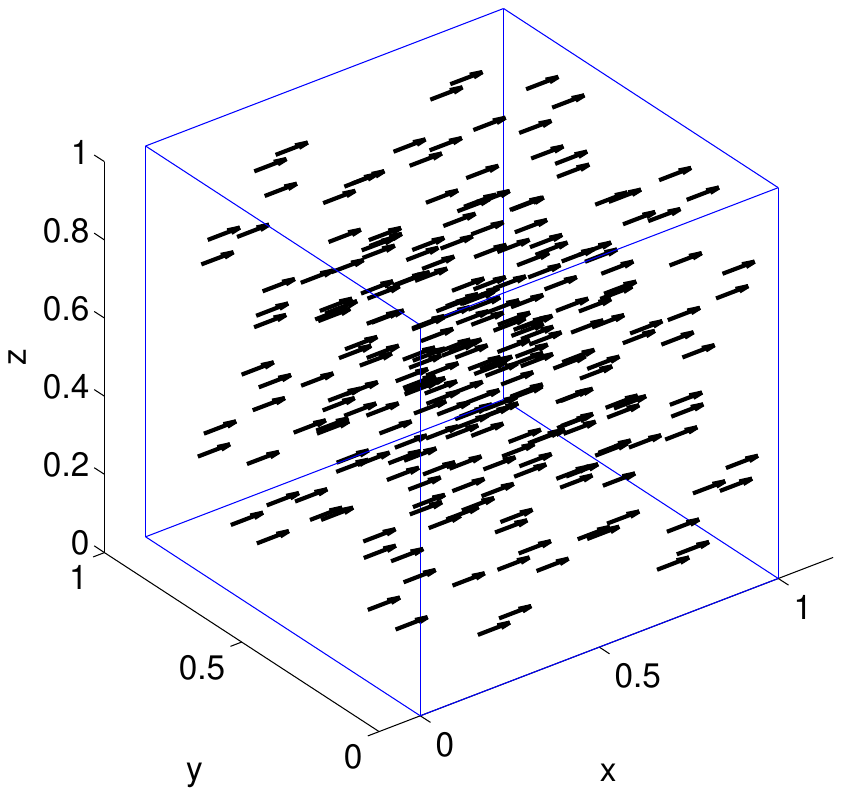}}
  \subfigure{
  \includegraphics[scale=0.35,trim=1in 2in 1.5in 2.2in,clip]
  {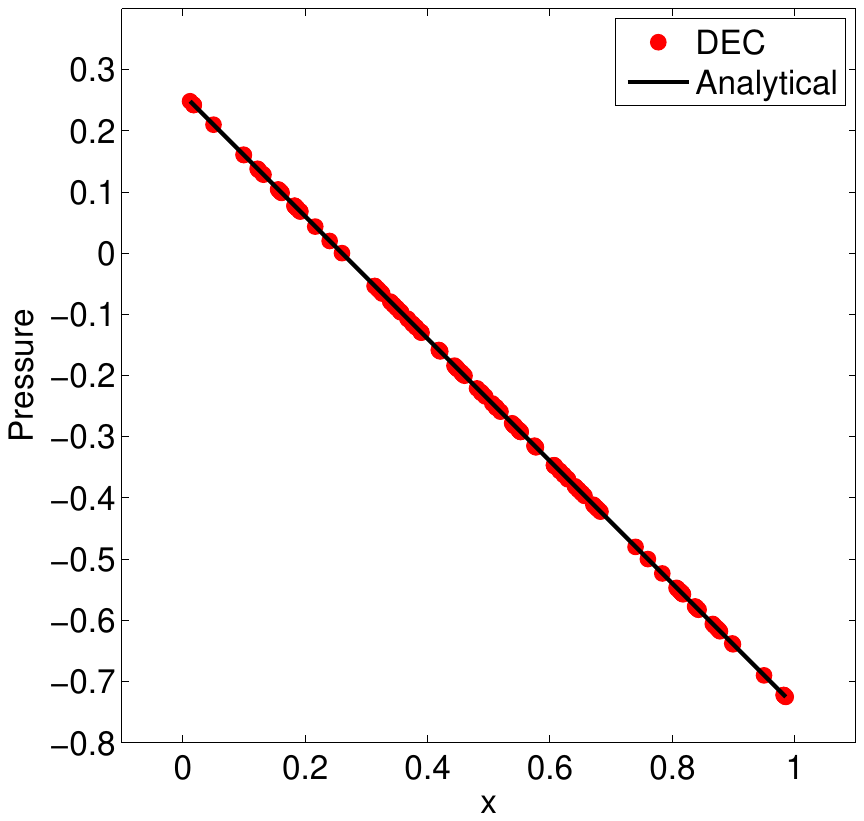}}
  \subfigure{
  \includegraphics[scale=0.35,trim=0.1in 2in 1in 1.1in,clip]
  {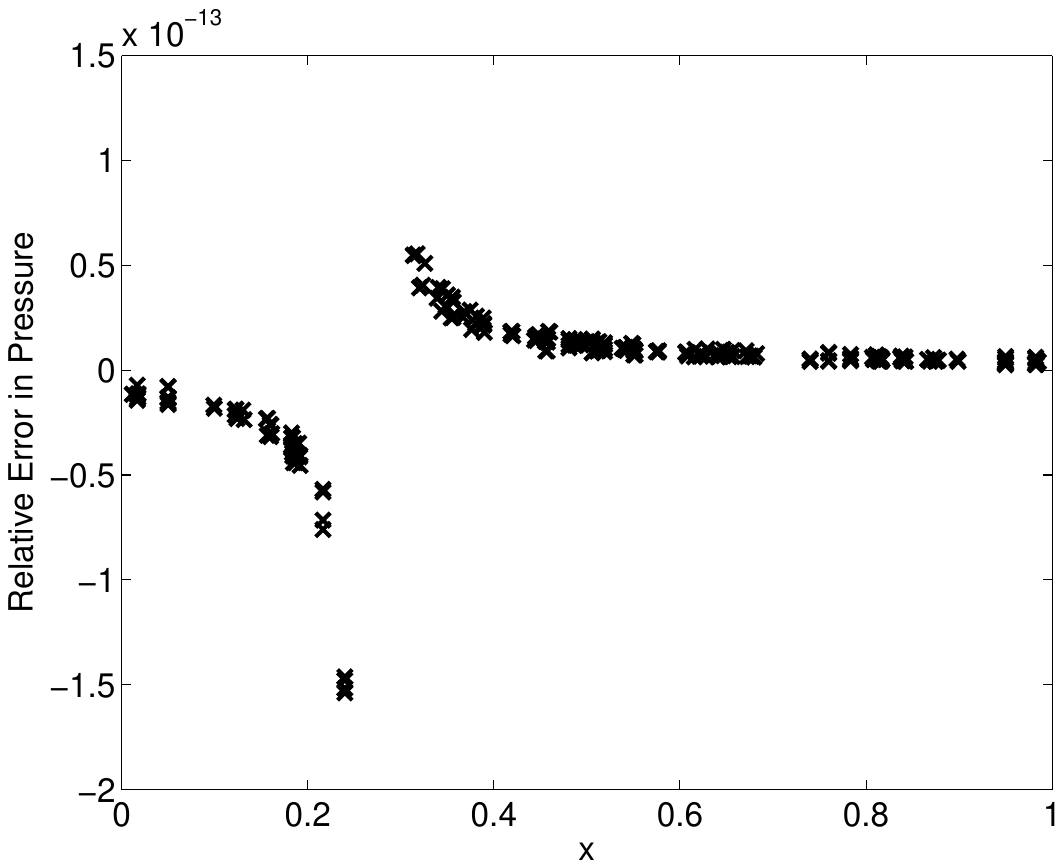}}
  \caption{Patch test on a cube. The mesh used for this cube has 244
    tetrahedra. For clarity, the tetrahedra in the
    cube have not been shown. The fluid velocity on the boundary is
    from negative $x$ to positive $x$ direction, with no $y$ or $z$
    components. The fluxes (2-cochains) across the internal faces are
    computed using our method and interpolated using the Whitney map.
    This is then sampled at the circumcenters, converted into a vector
    field and plotted as arrows. In the relative error plot for
    pressure, the data with 0 exact pressure has been removed.}
  \label{fig:cube244cnstvel}
\end{figure}

One of the focuses of this paper is on structure preserving numerical
methods, generated by a systematic application of DEC. The proposed
method has both local and global mass balance properties. In addition,
we show that the proposed method can also exactly represent linear
variation of pressure within a domain. We also highlight the need for
careful choice of locations for the pressure to get better numerical
solutions. DEC naturally provides the locations of pressure that
conserve local and global mass balance, and also exactly represent
linear variation of pressure. In DEC, for each element, pressure is
located at circumcenter. The fact that each $(n-1)$-simplex and its
circumcentric dual edge are orthogonal is important here. One common
choice used in the literature as location points for pressure are
barycenters. In Figure~\ref{fig:bary_vs_circum} we show that in DEC,
the choice of barycenter cannot exactly represent linear variation of
pressure (along with local and global mass balance) whereas location
of pressure at circumcenter can. However, we must point out that there
are other numerical schemes which are exact for linear solutions when
the barycenter is used. Thus the good behavior of circumcentric
location for pressure is not a general feature of numerical methods
but limited to DEC.

\subsection{Convergence tests in 2D} \label{subsec:cnvrgnc}
\begin{figure}[p]
  \centering
  \includegraphics[scale=0.6,trim=0in 0in 0in 0in, clip]
  {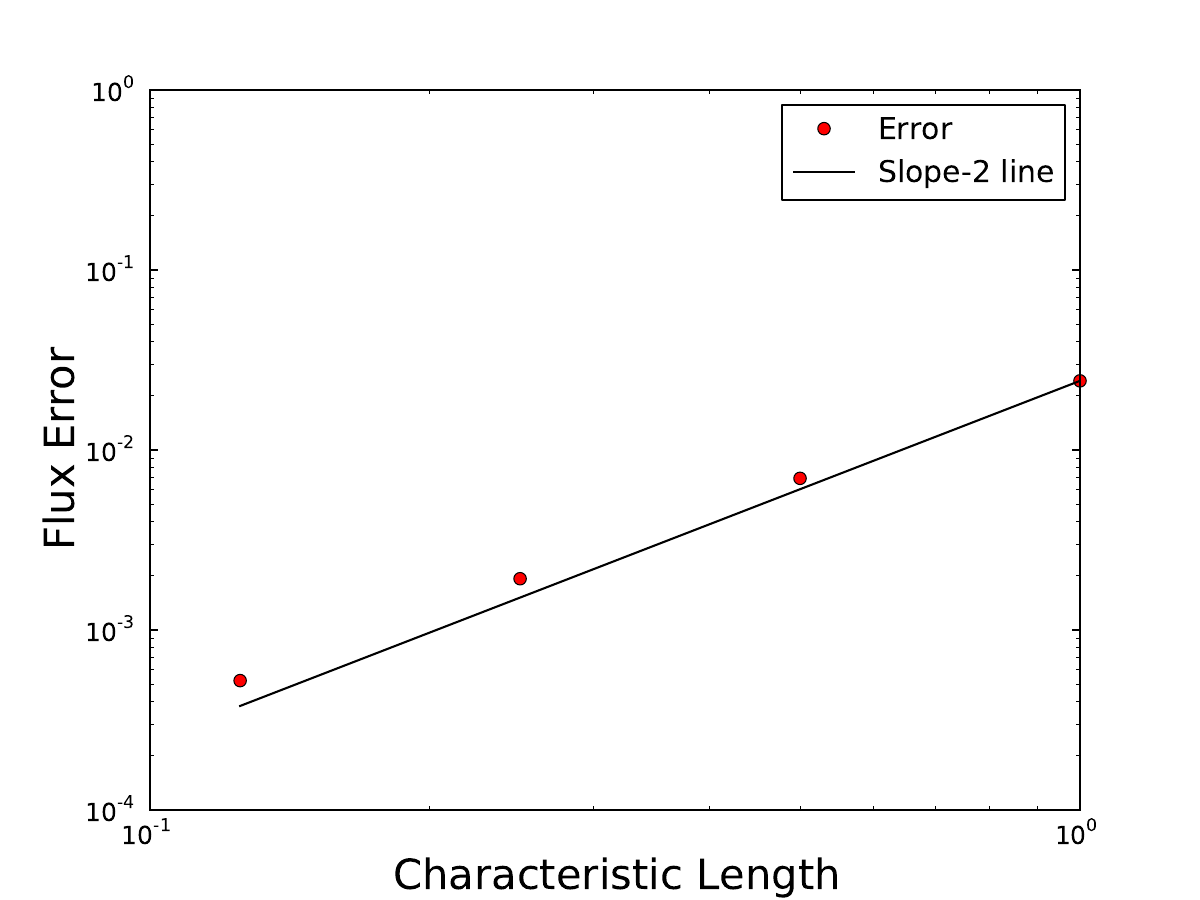}
  \includegraphics[scale=0.6,trim=0in 0in 0in 0in, clip]
  {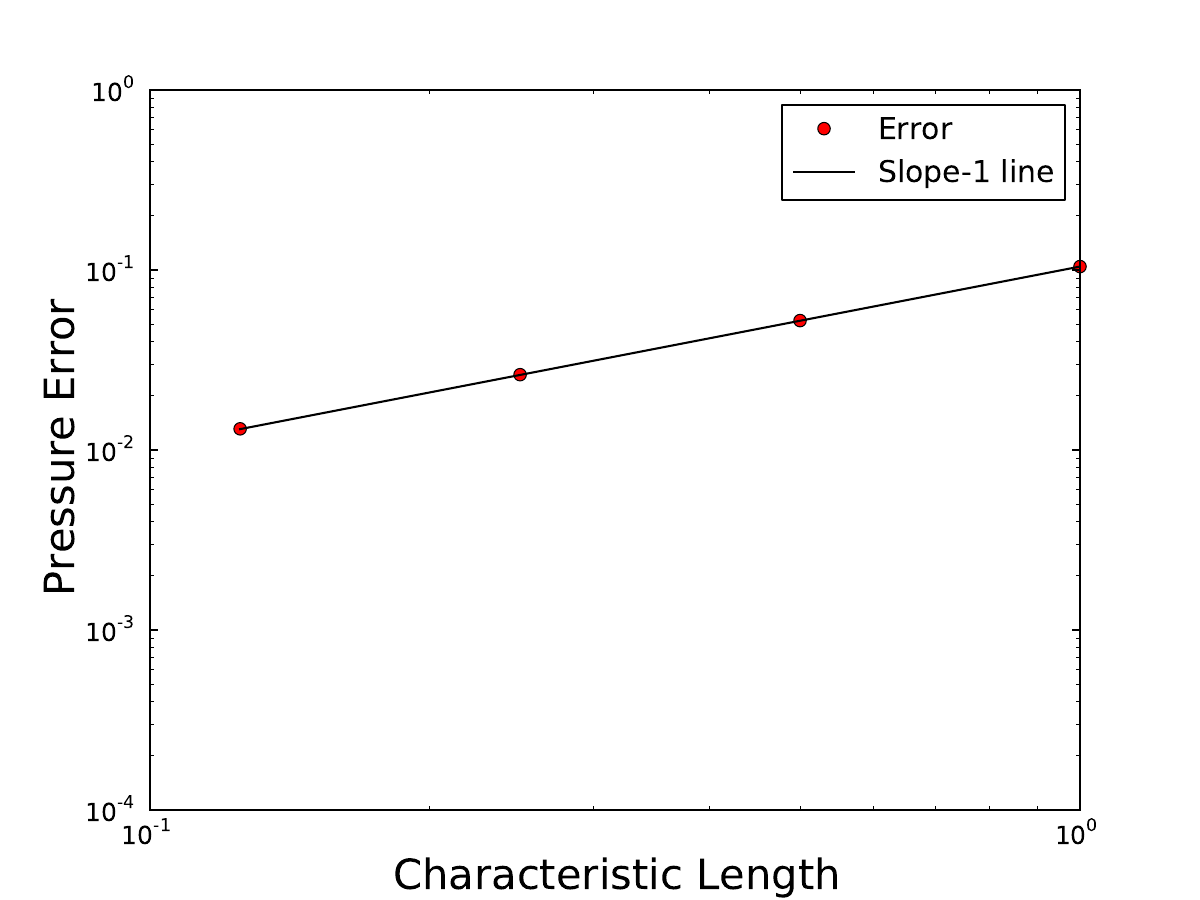}
  \caption{Numerical test of convergence for problem of
    Figure~\ref{fig:sqr336coscos}, starting with a mesh with 186
    triangles and subdividing it 3 times to obtain a sequence of 4
    meshes. A loglog plot of the norm of the absolute error in flux is
    shown on top. The bottom plot shows the absolute error in the
    pressure. See the text for how the errors are measured. The
    average slope of the flux error plot is about 1.9 and that for the
    pressure is about 1.03 as would be expected.}
\label{fig:flx_cnvrgnc}
\end{figure}

We used the setup of Section~\ref{subsec:general}
(Figure~\ref{fig:sqr336coscos}) for numerical test of convergence in
2D since the analytical solution is known. The resulting error plot as
a function of mesh length parameter is shown in
Figure~\ref{fig:flx_cnvrgnc}. The initial mesh used was a square with
186 triangles. This was subdivided 3 times to obtain a sequence of 3
finer meshes. The subdivision used was Loop
subdivision~\cite{Loop1987} in which each triangle is subdivided into
4 triangles similar to the original one and congruent to each other.

Since the flux is a primal 1-cochain it is associated with the primal
edges. Since DEC is like a finite difference or finite volume method,
it is appropriate to measure the flux error at the edge and then
integrate the error over the mesh. We compute the true flux through
the edge by using numerical quadrature. The difference between the
true and computed fluxes gives one number for each internal edge. We
interpolate the squares of these numbers from the edges to the
triangles using Whitney forms as explained in
Section~\ref{subsec:fluxviz}. The integral of this quantity can be
computed exactly. The square root of this is what is plotted in
Figure~\ref{fig:flx_cnvrgnc}. Likewise, pressure which is a dual
0-cochain is taken to be constant over each associated primal
2-simplex. We compute the error as the normed difference from the true
pressure integrated over the entire mesh. We determine this via
quadrature over each 2-simplex and computing the 2-norm of the
resulting vector of errors over 2-simplices. This is the quantity
plotted in Figure~\ref{fig:flx_cnvrgnc}. The numerical order of
convergence for flux is about 1.9 and that for pressure is about 1.04.


\subsection{Flow on a surface} \label{subsec:surface}

In recent years there has been considerable interest in surface finite
element methods. See for instance~\cite{DeDz2007,Demlow2009} and the
references therein. For surface methods, most theoretical work and all
the numerical experiments in the literature are for scalar elliptic
equations. In~\cite{HoSt2010} the authors develop an abstract
framework for analyzing mixed methods for equations involving
Laplace-deRham operators in arbitrary dimensions. However, there are
no specific finite elements proposed and no experimental results
shown. Recently, the first author of this paper and Kalyanaraman used
the method proposed in this paper (and its generalization to finite
element exterior calculus) for experiments on Darcy flow on a
surface~\cite{HiKa2011}. Since Darcy flow is equivalent to Poisson's
equation in mixed form, the results in~\cite{HiKa2011} are the first
experimental results in this setting. We reproduce some of these
experiments to show how easily our formulation generalizes to
surfaces. Indeed, no change in the programs is required in moving from
planar meshes to surface meshes.

\begin{figure}[p]
  \centering
  \includegraphics[scale=0.7, trim=3in 1in 3in 2in, clip]
  {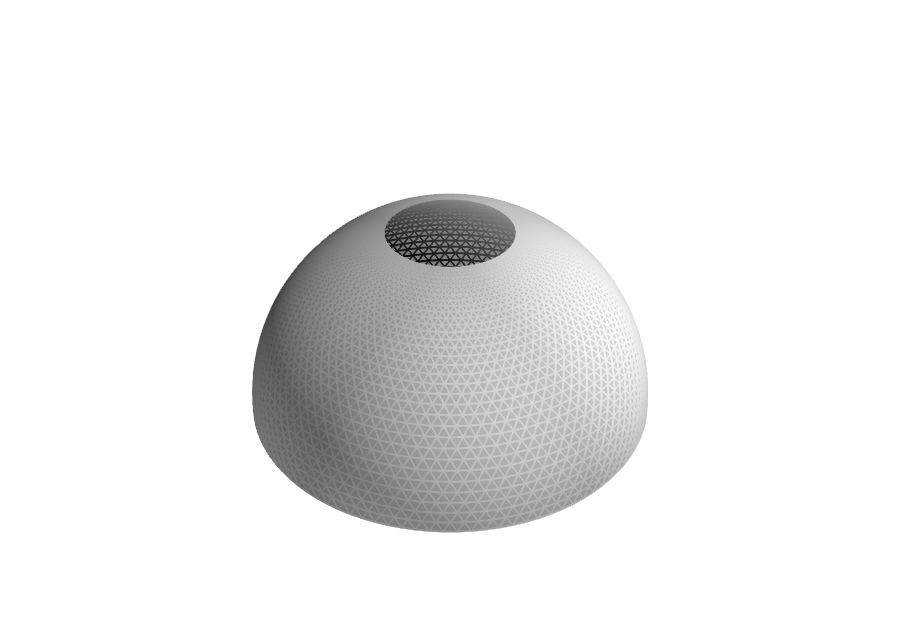}
  \caption{An example of a well-centered triangulation of an annular
    hemisphere. The mesh consists of 6600 triangles.}
  \label{fig:hmsphrhlmsh}
\end{figure}

\begin{figure}[!htp]
  \centering
  \begin{tabular}{c}
    \includegraphics[scale=0.6,trim=0in 0in 0in 0in, clip]
    {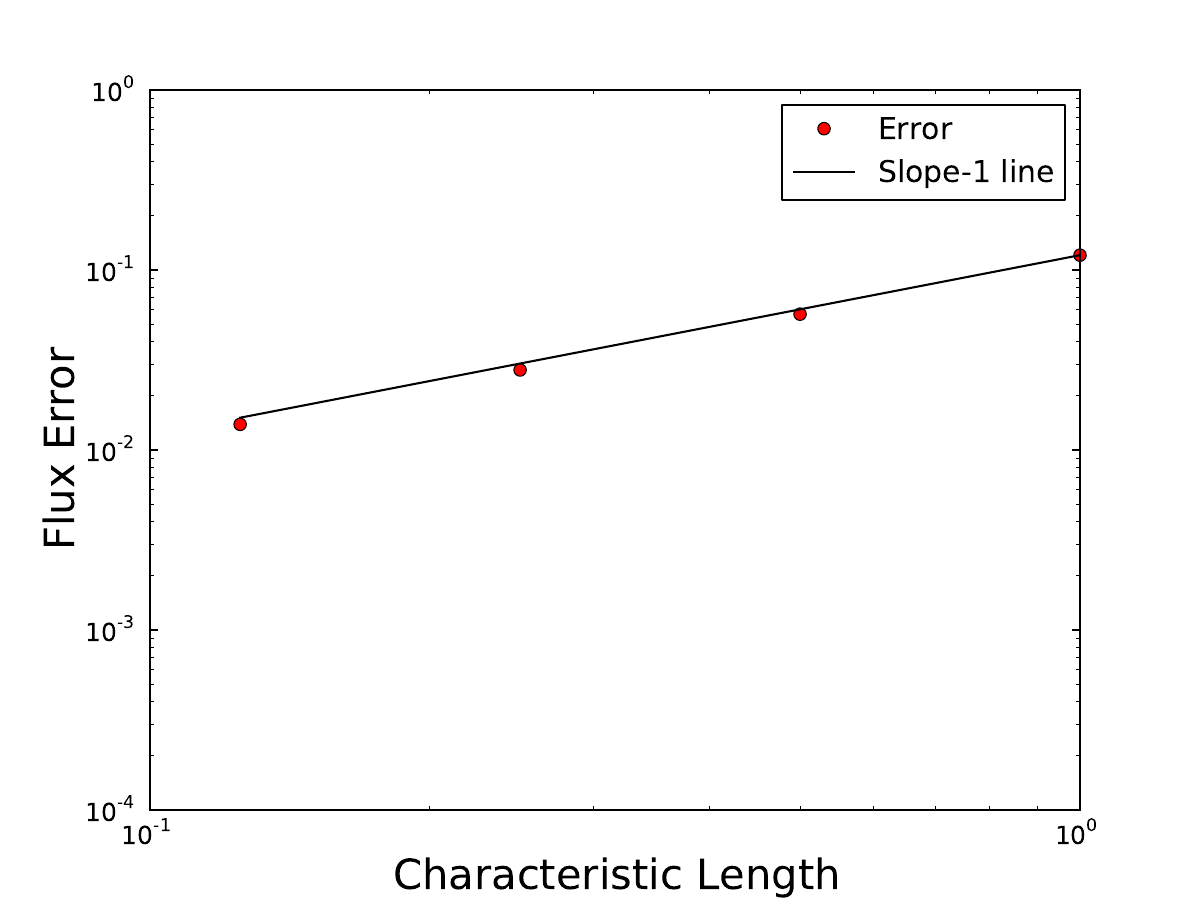} \\
    \includegraphics[scale=0.6,trim=0in 0in 0in 0in, clip]
    {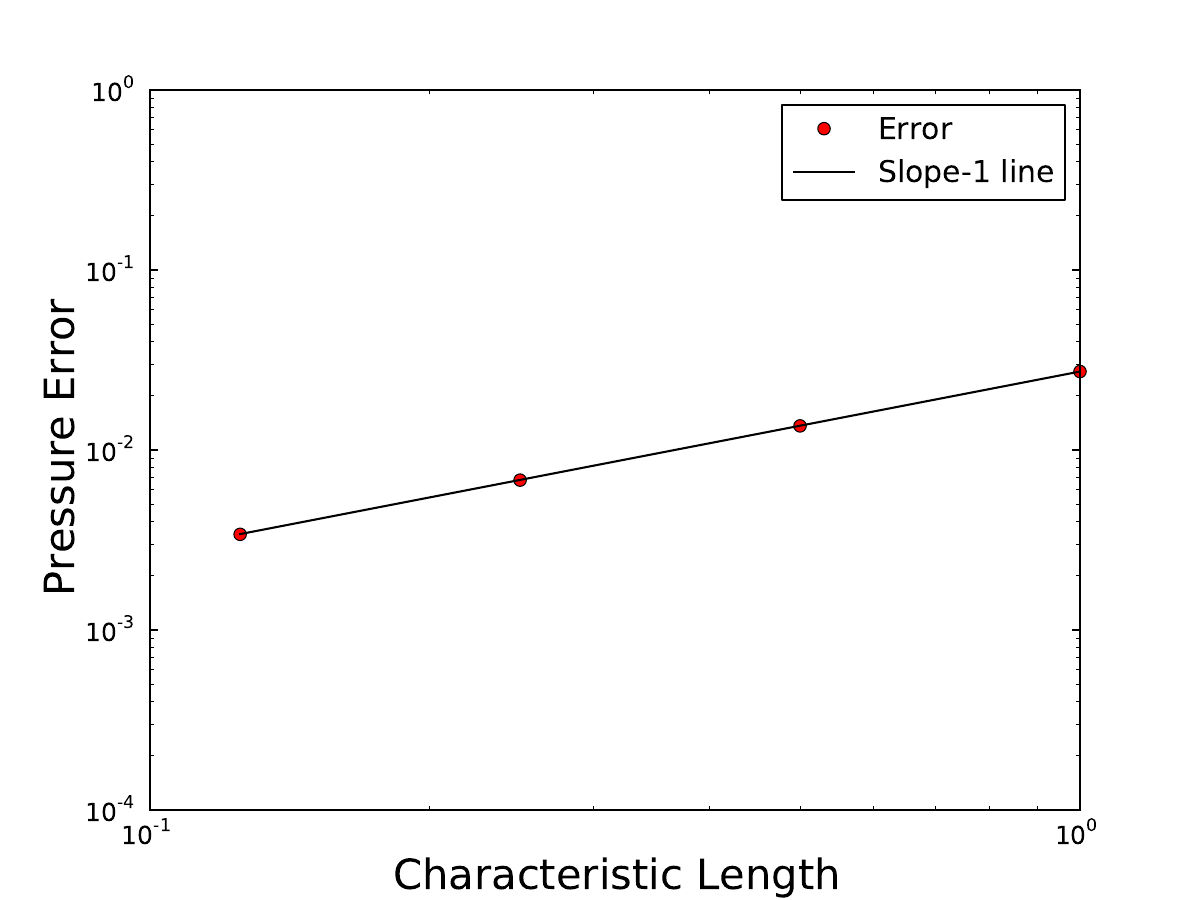}
  \end{tabular}
  \caption{Numerical test of convergence for problem of
    Figure~\ref{fig:hmsphrhlmsh}, starting with a mesh with 960
    triangles and subdividing it 3 times to obtain a sequence of 4
    meshes. A loglog plot of the norm of the absolute error in flux is
    shown on top. The bottom plot shows the absolute error in the
    pressure. See the text for how the errors are measured. The
    average slope of the flux error plot is about 1.04 and that for
    the pressure is about 1.0 as would be expected.}
  \label{fig:hmsphrhlDECcnvrgnc}
\end{figure}

Figure~\ref{fig:hmsphrhlmsh} shows an example of the family of
triangle mesh surfaces on which the experiments were done. This is a
piecewise linear approximation of a hemisphere with a hole cut out
near the north pole. The inflow is from the top circular boundary and
the outflow is from the bottom boundary and both are tangential to the
surface. The inflow is constant in magnitude around the circular
boundary and by symmetry so is the outflow, which can be computed
analytically. Figure~\ref{fig:hmsphrhlDECcnvrgnc} shows the numerical
evidence for convergence of our method in this setting. The
computation of the flux and pressure errors was done in a way similar
to the one described in Section~\ref{subsec:cnvrgnc}. The edges of the
triangulation were projected radially to the hemisphere. The flux
through the resulting curvilinear edges was computed by high order
adaptive numerical quadrature to get a true value of a flux through
each edge. The error in flux is then integrated on the piecewise
linear mesh as in Section~\ref{subsec:cnvrgnc}. The true pressure on
the sphere is projected radially to the triangles of the mesh and
compared with the constant pressures in those. The pressure error is
computed using a fourth order quadrature on each triangle.

\section{Conclusions and Future Work}
\label{sec:conclusions}

Our goal in this paper has been to introduce a numerical method for
Darcy flow derived using discrete exterior calculus (DEC). We have
shown that this approach results in a unified derivation of methods
for both two- and three-dimensions and for surfaces embedded in three
dimensions. We have also demonstrated the numerical performance of
these methods. For example, the method is shown to pass several patch
and other standard test problems in both two and three dimensions. We
also showed numerical convergence of pressure as well as flux. Our DEC
based approach is intrinsic, in the sense that it involves quantities
and operations that do not depend on how the mesh is embedded in
$\R^3$. As a result, it is easy to use the method for Darcy flow on a
curved surface using a discretization of such a surface as was
demonstrated. The formulation is such that generalization to finite
element exterior calculus involves the use of Whitney Hodge star
instead of the primal-dual DEC Hodge star used in this
paper. Experiments in this direction are an avenue for future work,
some of which has been explored in~\cite{HiKa2011}. Our method also
generalizes in another way, to a dual formulation, and this has been
done in~\cite{GiBa2010}.

Even in the case of domains that are open subsets of $\R^2$ or $\R^3$
the required computations such as circumcenter calculation, and
computation of volumes, areas and lengths, can be built from linear
algebra operations which simplifies the implementation as demonstrated
in \cite{BeHi2011}.

The philosophy of DEC is to discretize the operators and objects in
such a way that the their properties in smooth calculus have discrete
analogs. As a result, for example, in our method mass balance is
satisfied both locally (that is, in each element) and globally because
of the way the fluxes are represented and because Stokes' theorem is
true by definition in DEC.

A pedagogical advantage of the DEC approach is that once the language
of exterior calculus has been mastered and the discretizations
understood, the translation from a PDE to its corresponding discrete
equation, and from there to the matrix form are trivial steps. Compare
for example the smooth equation~\eqref{eq:dfdrcy}, the corresponding
discrete equation~\eqref{eq:darcy} and the first block row of the
matrix equation~\eqref{eq:mtxwthIfnl}. There is also a clean
separation in the matrix form equation~\eqref{eq:mtxwthstr} between
the parts depending only on the mesh connectivity (the top right and
bottom left blocks) and that depending on how the physical variable
are related (the top left part) which usually involves constitutive
relationships as well. These depend on material measurements, and
hence the inaccuracies of measurement of material properties does not
corrupt that part of the matrix that is topological. In this aspect
our method shares this good property of some mixed finite element
method formulations.

In addition to applying DEC to Darcy flow, in this paper we have also
developed a discretization of Hodge star for non-homogeneous
medium. We used this in examples in which the permeability varies
across the mesh, either along the flow or transversal to the
flow. This discretization of a spatially varying Hodge star should be
useful in other PDEs as well.

There are many avenues for future research even within Darcy flow.
Flow on curved surfaces has already been mentioned before. Another
direction for further work is the proper DEC discretization of
anisotropic permeability. Finally the analytical treatment of
convergence and stability properties remains to be done. For this, it
might be useful to formulate this method as a finite element method.

\section*{Acknowledgments}
The research of ANH and JHC was supported by the National Science
Foundation with an NSF CAREER Award (Grant No. DMS-0645604). The
research of KBN was supported by the Texas Engineering Experiment
Station. ANH would like to thank his student Evan VanderZee for help
in creating some of the meshes in this paper. ANH would also like to
thank his student Kaushik Kalyanaraman for help with the convergence
tests in 2D and for the experiments on surface meshes. The opinions
expressed in this paper are those of the authors and do not
necessarily reflect that of the sponsors.

\bibliographystyle{acmurldoi}
\bibliography{hirani}

\begin{thebibliography}{10}

\bibitem{AbMaRa1988}
{\sc Abraham, R., Marsden, J.~E., and Ratiu, T.}
\newblock {\em Manifolds, Tensor Analysis, and Applications}, second~ed.
\newblock Springer--Verlag, New York, 1988.

\bibitem{AcBeCo2003}
{\sc Achdou, Y., Bernardi, C., and Coquela, F.}
\newblock A priori and a posteriori analysis of finite volume discretizations
  of {D}arcy's equations.
\newblock {\em Numerische Mathematik 96\/} (2003), 17--42.

\bibitem{ArBoLeNiSh2006}
{\sc Arnold, D.~N., Bochev, P.~B., Lehoucq, R.~B., Nicolaides, R.~A., and
  Shashkov, M.}, Eds.
\newblock {\em Compatible Spatial Discretizations}, vol.~142 of {\em The IMA
  Volumes in Mathematics and its Applications}.
\newblock Springer New York, 2006.

\bibitem{ArFaWi2006}
{\sc Arnold, D.~N., Falk, R.~S., and Winther, R.}
\newblock Finite element exterior calculus, homological techniques, and
  applications.
\newblock In {\em Acta Numerica}, A.~Iserles, Ed., vol.~15. Cambridge
  University Press, 2006, pp.~1--155.

\bibitem{ArFaWi2010}
{\sc Arnold, D.~N., Falk, R.~S., and Winther, R.}
\newblock Finite element exterior calculus: from {H}odge theory to numerical
  stability.
\newblock {\em Bull. Amer. Math. Soc. (N.S.) 47}, 2 (2010), 281--354.
\newblock \href {http://dx.doi.org/10.1090/S0273-0979-10-01278-4} {\path{doi:
  10.1090/S0273-0979-10-01278-4}}.

\bibitem{Arnold1989}
{\sc Arnold, V.~I.}
\newblock {\em Mathematical Methods of Classical Mechanics}, second~ed.
\newblock Springer--Verlag, New York, 1989.
\newblock Translated from the Russian by K. Vogtmann and A. Weinstein.

\bibitem{BaNa1997}
{\sc Babuska, I., and Narasimhan, R.}
\newblock The {B}abuska-{B}rezzi condition and the patch test: an example.
\newblock {\em Computer Methods in Applied Mechanics and Engineering 140}, 1-2
  (January 1997), 183--199.

\bibitem{BeHi2011}
{\sc Bell, N., and Hirani, A.~N.}
\newblock {PyDEC}: Algorithms and software for {D}iscretization of {E}xterior
  {C}alculus, March 2011.
\newblock Available as e-print on arxiv.org.
\newblock URL \url{http://arxiv.org/abs/1103.3076}, \href
  {http://arxiv.org/abs/1103.3076v1} {\path{arXiv:1103.3076v1}}.

\bibitem{BeGoLi2005}
{\sc Benzi, M., Golub, G.~H., and Liesen, J.}
\newblock Numerical solution of saddle point problems.
\newblock {\em Acta Numerica 14\/} (2005), 1--137.

\bibitem{BoDo2006}
{\sc Bochev, P., and Dohrmann, C.}
\newblock A computational study of stabilized, low-order {$C^0$} finite element
  approximations of {D}arcy equations.
\newblock {\em Computational Mechanics 38}, 4-5 (2006), 323--333.

\bibitem{BoHy2006}
{\sc Bochev, P.~B., and Hyman, J.~M.}
\newblock Principles of mimetic discretizations of differential operators.
\newblock In {\em Compatible Spatial Discretizations}, D.~N. Arnold, P.~B.
  Bochev, R.~B. Lehoucq, R.~A. Nicolaides, and M.~Shashkov, Eds., vol.~142 of
  {\em The IMA Volumes in Mathematics and its Applications}. Springer, Berlin,
  2006, pp.~89--119.

\bibitem{Bossavit1988b}
{\sc Bossavit, A.}
\newblock Mixed finite elements and the complex of {W}hitney forms.
\newblock In {\em The Mathematics of Finite Elements and Applications {VI}},
  J.~Whiteman, Ed. Academic Press, 1988, pp.~137--144.

\bibitem{Bossavit1988}
{\sc Bossavit, A.}
\newblock A rationale for ``edge elements" in 3-{D} fields computations.
\newblock {\em IEEE Trans. Mag. 24}, 1 (January 1988), 74--79.

\bibitem{Bossavit1988a}
{\sc Bossavit, A.}
\newblock Whitney forms : {A} class of finite elements for three-dimensional
  computations in electromagnetism.
\newblock {\em IEE Proceedings 135, Part A}, 8 (November 1988), 493--500.

\bibitem{Bossavit1998}
{\sc Bossavit, A.}
\newblock {\em Computational Electromagnetism : Variational Formulations,
  Complementarity, Edge Elements}.
\newblock Academic Press, 1998.

\bibitem{Bossavit1998e}
{\sc Bossavit, A.}
\newblock On the geometry of electromagnetism (4): {M}axwell's house.
\newblock {\em Journal of the Japan Society of Applied Electromagnetics 6}, 4
  (1998), 318--326.

\bibitem{Braess2007}
{\sc Braess, D.}
\newblock {\em Finite Elements: Theory, Fast Solvers, and Applications in Solid
  Mechanics}, third~ed.
\newblock Cambridge University Press, Cambridge, 2007.
\newblock Translated from the 1992 German edition by Larry L. Schumaker.

\bibitem{BrDoMa1985}
{\sc Brezzi, F., Douglas, J., and Marini, L.~D.}
\newblock Two families of mixed finite elements for second order elliptic
  problems.
\newblock {\em Numerische Mathematik 47}, 2 (June 1985), 217--235.

\bibitem{BrFo1991}
{\sc Brezzi, F., and Fortin, M.}
\newblock {\em {Mixed and hybrid finite element methods, volume 15 of Springer
  series in computational mathematics}}.
\newblock Springer-Verlag, New York, 1991.

\bibitem{BrLiSh2005}
{\sc Brezzi, F., Lipnikov, K., and Shashkov, M.}
\newblock Convergence of the mimetic finite difference method for diffusion
  problems on polyhedral meshes.
\newblock {\em SIAM Journal on Numerical Analysis 43}, 5 (2005), 1872--1896.
\newblock URL \url{http://link.aip.org/link/?SNA/43/1872/1}, \href
  {http://dx.doi.org/10.1137/040613950} {\path{doi: 10.1137/040613950}}.

\bibitem{BrLiShSi2007}
{\sc Brezzi, F., Lipnikov, K., Shashkov, M., and Simoncini, V.}
\newblock A new discretization methodology for diffusion problems on
  generalized polyhedral meshes.
\newblock {\em Computer Methods in Applied Mechanics and Engineering 196},
  37-40 (2007), 3682 -- 3692.
\newblock Special Issue Honoring the 80th Birthday of Professor Ivo Babuska.
\newblock URL
  \url{http://www.sciencedirect.com/science/article/B6V29-4N8M857-1/2/a6cacc4af91c15ee1a16ac542e906dde},
  \href {http://dx.doi.org/DOI: 10.1016/j.cma.2006.10.028} {\path{doi: DOI:
  10.1016/j.cma.2006.10.028}}.

\bibitem{Chen2005}
{\sc Chen, Z.}
\newblock {\em Finite element methods and their applications}.
\newblock Scientific Computation. Springer-Verlag, Berlin, 2005.

\bibitem{ChVa1999}
{\sc Chou, S., and Vassilevski, P.}
\newblock A general mixed covolume framework for constructing conservative
  schemes for elliptic problems.
\newblock {\em Mathematics of Computation 68}, 227 (1999), 991--1011.

\bibitem{ChKwVa1998}
{\sc Chou, S.-H., Kwak, D.~Y., and Vassilevski, P.~S.}
\newblock Mixed covolume methods for elliptic problems on triangular grids.
\newblock {\em SIAM Journal on Numerical Analysis 35}, 5 (1998), 1850--1861.

\bibitem{Demlow2009}
{\sc Demlow, A.}
\newblock Higher-order finite element methods and pointwise error estimates for
  elliptic problems on surfaces.
\newblock {\em SIAM Journal on Numerical Analysis 47}, 2 (2009), 805--827.
\newblock URL \url{http://link.aip.org/link/?SNA/47/805/1}, \href
  {http://dx.doi.org/10.1137/070708135} {\path{doi: 10.1137/070708135}}.

\bibitem{DeDz2007}
{\sc Demlow, A., and Dziuk, G.}
\newblock An adaptive finite element method for the {L}aplace--{B}eltrami
  operator on implicitly defined surfaces.
\newblock {\em SIAM Journal on Numerical Analysis 45}, 1 (2007), 421--442.
\newblock URL \url{http://link.aip.org/link/?SNA/45/421/1}, \href
  {http://dx.doi.org/10.1137/050642873} {\path{doi: 10.1137/050642873}}.

\bibitem{DeHiLeMa2005}
{\sc Desbrun, M., Hirani, A.~N., Leok, M., and Marsden, J.~E.}
\newblock Discrete exterior calculus, August 2005.
\newblock Available as e-print on arxiv.org.
\newblock \href {http://arxiv.org/abs/math.DG/0508341}
  {\path{arXiv:math.DG/0508341}}.

\bibitem{DeKaTo2008}
{\sc Desbrun, M., Kanso, E., and Tong, Y.}
\newblock Discrete differential forms for computational modeling.
\newblock In {\em Discrete Differential Geometry}, A.~I. Bobenko, J.~M.
  Sullivan, P.~Schr\"oder, and G.~M. Ziegler, Eds., vol.~38 of {\em Oberwolfach
  Seminars}. Birkh\"auser Basel, 2008, pp.~287--324.
\newblock \href {http://dx.doi.org/10.1007/978-3-7643-8621-4_16} {\path{doi:
  10.1007/978-3-7643-8621-4_16}}.

\bibitem{Dodziuk1976}
{\sc Dodziuk, J.}
\newblock Finite-difference approach to the {H}odge theory of harmonic forms.
\newblock {\em Amer. J. Math. 98}, 1 (1976), 79--104.

\bibitem{DuGuJu2003a}
{\sc Du, Q., Gunzburger, M.~D., and Ju, L.}
\newblock Voronoi-based finite volume methods, optimal {V}oronoi meshes, and
  {PDE}s on the sphere.
\newblock {\em Computer Methods in Applied Mechanics and Engineering 192},
  35-36 (August 2003), 3933--3957.
\newblock URL \url{http://dx.doi.org/10.1016/S0045-7825(03)00394-3}, \href
  {http://dx.doi.org/10.1016/S0045-7825(03)00394-3} {\path{doi:
  10.1016/S0045-7825(03)00394-3}}.

\bibitem{DuJu2005}
{\sc Du, Q., and Ju, L.}
\newblock Finite volume methods on spheres and spherical centroidal voronoi
  meshes.
\newblock {\em SIAM Journal on Numerical Analysis 43}, 4 (2005), 1673--1692.
\newblock URL \url{http://link.aip.org/link/?SNA/43/1673/1}, \href
  {http://dx.doi.org/10.1137/S0036142903425410} {\path{doi:
  10.1137/S0036142903425410}}.

\bibitem{Dziuk1988}
{\sc Dziuk, G.}
\newblock Finite elements for the {B}eltrami operator on arbitrary surfaces.
\newblock In {\em Partial Differential Equations and Calculus of Variations},
  vol.~Volume 1357/1988 of {\em Lecture Notes in Mathematics}. Springer Berlin
  Heidelberg, 1988, pp.~142--155.
\newblock URL \url{http://dx.doi.org/10.1007/BFb0082865}, \href
  {http://dx.doi.org/10.1007/BFb0082865} {\path{doi: 10.1007/BFb0082865}}.

\bibitem{ElToKaScDe2007}
{\sc Elcott, S., Tong, Y., Kanso, E., Schr\"oder, P., and Desbrun, M.}
\newblock Stable, circulation-preserving, simplicial fluids.
\newblock {\em ACM Transactions on Graphics 26}, 1 (2007), 4.

\bibitem{Frankel2004}
{\sc Frankel, T.}
\newblock {\em The Geometry of Physics}, second~ed.
\newblock Cambridge University Press, Cambridge, 2004.
\newblock An introduction.

\bibitem{GiBa2010}
{\sc Gillette, A., and Bajaj, C.}
\newblock A generalization for stable mixed finite elements.
\newblock In {\em SPM '10: Proceedings of the 14th ACM Symposium on Solid and
  Physical Modeling\/} (New York, NY, USA, 2010), ACM, pp.~41--50.
\newblock \href {http://dx.doi.org/10.1145/1839778.1839785} {\path{doi:
  10.1145/1839778.1839785}}.

\bibitem{HiPoWa2006}
{\sc Hildebrandt, K., Polthier, K., and Wardetzky, M.}
\newblock On the convergence of metric and geometric properties of polyhedral
  surfaces.
\newblock {\em Geom. Dedicata 123\/} (2006), 89--112.
\newblock URL \url{http://dx.doi.org/10.1007/s10711-006-9109-5}, \href
  {http://dx.doi.org/10.1007/s10711-006-9109-5} {\path{doi:
  10.1007/s10711-006-9109-5}}.

\bibitem{Hirani2003}
{\sc Hirani, A.~N.}
\newblock {\em Discrete Exterior Calculus}.
\newblock PhD thesis, California Institute of Technology, May 2003.
\newblock URL \url{http://resolver.caltech.edu/CaltechETD:etd-05202003-095403}.

\bibitem{HiKa2011}
{\sc Hirani, A.~N., and Kalyanaraman, K.}
\newblock Numerical experiments for {D}arcy flow on a surface using mixed
  exterior calculus methods, March 2011.
\newblock Available as e-print on arxiv.org.
\newblock \href {http://arxiv.org/abs/1103.4865} {\path{arXiv:1103.4865}}.

\bibitem{HoSt2010}
{\sc Holst, M., and Stern, A.}
\newblock Geometric variational crimes: Hilbert complexes, finite element
  exterior calculus, and problems on hypersurfaces.
\newblock {\em arXiv:1005.4455v1 [math.NA] Online\/} (May 2010).
\newblock URL \url{http://arxiv.org/abs/1005.4455}, \href
  {http://arxiv.org/abs/1005.4455} {\path{arXiv:1005.4455}}.

\bibitem{HuMaWa2006}
{\sc Hughes, T. J.~R., Masud, A., and Wan, J.}
\newblock A stabilized mixed discontinuous {G}alerkin method for {D}arcy flow.
\newblock {\em Computer Methods in Applied Mechanics and Engineering 195\/}
  (2006), 3347--3381.

\bibitem{HySh1997b}
{\sc Hyman, J.~M., and Shashkov, M.}
\newblock Adjoint operators for the natural discretizations of the divergence,
  gradient and curl on logically rectangular grids.
\newblock {\em Appl. Numer. Math. 25}, 4 (1997), 413--442.
\newblock URL \url{http://dx.doi.org/10.1016/S0168-9274(97)00097-4}, \href
  {http://dx.doi.org/10.1016/S0168-9274(97)00097-4} {\path{doi:
  10.1016/S0168-9274(97)00097-4}}.

\bibitem{HySh1997a}
{\sc Hyman, J.~M., and Shashkov, M.}
\newblock Natural discretizations for the divergence, gradient, and curl on
  logically rectangular grids.
\newblock {\em Comput. Math. Appl. 33}, 4 (1997), 81--104.
\newblock URL \url{http://dx.doi.org/10.1016/S0898-1221(97)00009-6}, \href
  {http://dx.doi.org/10.1016/S0898-1221(97)00009-6} {\path{doi:
  10.1016/S0898-1221(97)00009-6}}.

\bibitem{HySh1999}
{\sc Hyman, J.~M., and Shashkov, M.}
\newblock The orthogonal decomposition theorems for mimetic finite difference
  methods.
\newblock {\em SIAM J. Numer. Anal. 36}, 3 (1999), 788--818.
\newblock URL \url{http://dx.doi.org/10.1137/S0036142996314044}, \href
  {http://dx.doi.org/10.1137/S0036142996314044} {\path{doi:
  10.1137/S0036142996314044}}.

\bibitem{IrLo1983}
{\sc Irons, B., and Loikkanen, M.}
\newblock An engineers' defence of the patch test.
\newblock {\em International Journal for Numerical Methods in Engineering 19},
  9 (1983), 1391--1401.

\bibitem{JuDu2009}
{\sc Ju, L., and Du, Q.}
\newblock A finite volume method on general surfaces and its error estimates.
\newblock {\em Journal of Mathematical Analysis and Applications 352}, 2
  (2009), 645--668.
\newblock \href {http://dx.doi.org/10.1016/j.jmaa.2008.11.022} {\path{doi:
  10.1016/j.jmaa.2008.11.022}}.

\bibitem{LiShSv2006}
{\sc Lipnikov, K., Shashkov, M., and Svyatskiy, D.}
\newblock The mimetic finite difference discretization of diffusion problem on
  unstructured polyhedral meshes.
\newblock {\em Journal of Computational Physics 211}, 2 (January 2006),
  473--491.

\bibitem{LiShSvVa2007}
{\sc Lipnikov, K., Shashkov, M., Svyatskiy, D., and Vassilevski, Y.}
\newblock Monotone finite volume schemes for diffusion equations on
  unstructured triangular and shape-regular polygonal meshes.
\newblock {\em Journal of Computational Physics 227}, 1 (492--512 2007).

\bibitem{Loop1987}
{\sc Loop, C.~T.}
\newblock Smooth subdivision surfaces based on triangles.
\newblock Master's thesis, University of Utah, Department of Mathematics,
  August 1987.
\newblock URL \url{http://research.microsoft.com/~cloop/thesis.pdf}.

\bibitem{MaHu2002}
{\sc Masud, A., and Hughes, T. J.~R.}
\newblock A stabilized mixed finite element method for {D}arcy flow.
\newblock {\em Computer Methods Applied Mechanics and Engineering 191\/}
  (2002), 4341--4370.

\bibitem{Mattiussi1997}
{\sc Mattiussi, C.}
\newblock An analysis of finite volume, finite element, and finite difference
  methods using some concepts from algebraic topology.
\newblock {\em Journal of Computational Physics 133}, 2 (1997), 289 -- 309.
\newblock URL
  \url{http://www.sciencedirect.com/science/article/B6WHY-45S92HD-57/2/7ea87a6234c2cf93e4e9c663dd359db2},
  \href {http://dx.doi.org/DOI: 10.1006/jcph.1997.5656} {\path{doi: DOI:
  10.1006/jcph.1997.5656}}.

\bibitem{Mattiussi2002}
{\sc Mattiussi, C.}
\newblock A reference discretization strategy for the numerical solution of
  physical field problems.
\newblock {\em Advances in Imaging and Electron Physics 121\/} (2002),
  143--279.

\bibitem{MuMcPaDuToKaMaDe2010}
{\sc Mullen, P., McKenzie, A., Pavlov, D., Durant, L., Tong, Y., Kanso, E.,
  Marsden, J.~E., and Desbrun, M.}
\newblock Discrete lie advection of differential forms.
\newblock {\em Foundations of Computational Mathematics In press\/} (2010).

\bibitem{Munkres1984}
{\sc Munkres, J.~R.}
\newblock {\em Elements of Algebraic Topology}.
\newblock Addison--Wesley Publishing Company, Menlo Park, 1984.

\bibitem{NaMaHj2008}
{\sc Nakshatrala, K.~B., Masud, A., and Hjelmstad, K.~D.}
\newblock On finite element formulations for nearly incompressible linear
  elasticity.
\newblock {\em Computational Mechanics 41\/} (2008), 547--561.

\bibitem{NaTuHjMa2006}
{\sc Nakshatrala, K.~B., Turner, D.~Z., Hjelmstad, K.~D., and Masud, A.}
\newblock A mixed stabilized finite element formulation for {D}arcy flow based
  on a multiscale decomposition of the solution.
\newblock {\em Computer Methods in Applied Mechanics and Engineering 195\/}
  (2006), 4036--4049.

\bibitem{Nedelec1980}
{\sc Nedelec, J.~C.}
\newblock Mixed finite elements in $\mathbb{R}^3$.
\newblock {\em Numerische Mathematik 35}, 3 (1980), 315--341.

\bibitem{Nicolaides1992}
{\sc Nicolaides, R.~A.}
\newblock Direct discretization of planar div-curl problems.
\newblock {\em SIAM Journal on Numerical Analysis 29}, 1 (1992), 32--56.
\newblock URL \url{http://www.jstor.org/stable/2158074}.

\bibitem{NiTr2006}
{\sc Nicolaides, R.~A., and Trapp, K.~A.}
\newblock Covolume discretization of differential forms.
\newblock In {\em Compatible Spatial Discretizations}, D.~N. Arnold, P.~B.
  Bochev, R.~B. Lehoucq, R.~A. Nicolaides, and M.~Shashkov, Eds., vol.~142 of
  {\em The IMA Volumes in Mathematics and its Applications}. Springer, New
  York, 2006, pp.~161--171.

\bibitem{PaMuToKaMaDe2009}
{\sc Pavlov, D., Mullen, P., Tong, Y., Kanso, E., Marsden, J.~E., and Desbrun,
  M.}
\newblock Structure-preserving discretization of incompressible fluids.
\newblock Available as e-print arXiv:0912.3989v2 [math.DS], December 20 2009.
\newblock URL \url{http://arxiv.org/abs/0912.3989v2}.

\bibitem{PeSu2007}
{\sc Perot, J., and Subramanian, V.}
\newblock Discrete calculus methods for diffusion.
\newblock {\em Journal of Computational Physics 224}, 1 (2007), 59--81.

\bibitem{RaTh1977}
{\sc Raviart, P.-A., and Thomas, J.~M.}
\newblock A mixed finite element method for 2nd order elliptic problems.
\newblock In {\em Mathematical aspects of finite element methods (Proc. Conf.,
  Consiglio Naz. delle Ricerche (C.N.R.), Rome, 1975)}, Lecture Notes in Math.,
  Vol. 606. Springer, Berlin, 1977, pp.~292--315.

\bibitem{RiJuGu2008}
{\sc Ringler, T., Ju, L., and Gunzburger, M.}
\newblock A multiresolution method for climate system modeling: application of
  spherical centroidal {V}oronoi tessellations.
\newblock {\em Ocean Dynamics 58}, 5-6 (2008), 475--498.
\newblock \href {http://dx.doi.org/10.1007/s10236-008-0157-2} {\path{doi:
  10.1007/s10236-008-0157-2}}.

\bibitem{Shewchuk1996}
{\sc Shewchuk, J.~R.}
\newblock Triangle: {E}ngineering a {2D} {Q}uality {M}esh {G}enerator and
  {D}elaunay {T}riangulator.
\newblock In {\em Applied Computational Geometry: Towards Geometric
  Engineering}, M.~C. Lin and D.~Manocha, Eds., vol.~1148 of {\em Lecture Notes
  in Computer Science}. Springer-Verlag, May 1996, pp.~203--222.
\newblock From the First ACM Workshop on Applied Computational Geometry.

\bibitem{Si2009}
{\sc Si, H.}
\newblock {TetGen} : {A} quality tetrahedral mesh generator and a 3{D}
  {D}elaunay triangulator [online].
\newblock 2009.
\newblock URL \url{http://tetgen.berlios.de/}.

\bibitem{VaHiGuRa2010}
{\sc VanderZee, E., Hirani, A.~N., Guoy, D., and Ramos, E.~A.}
\newblock Well-centered triangulation.
\newblock {\em SIAM Journal on Scientific Computing 31}, 6 (2010), 4497--4523.
\newblock \href {http://dx.doi.org/10.1137/090748214} {\path{doi:
  10.1137/090748214}}.

\bibitem{Wardetzky2008}
{\sc Wardetzky, M.}
\newblock Convergence of the cotangent formula: {A}n overview.
\newblock In {\em Discrete Differential Geometry}, A.~I. Bobenko, J.~M.
  Sullivan, P.~Schr\"oder, and G.~M. Ziegler, Eds., vol.~38 of {\em Oberwolfach
  Seminars}. Birkh\"auser Basel, 2008, pp.~275--286.
\newblock \href {http://dx.doi.org/10.1007/978-3-7643-8621-4_15} {\path{doi:
  10.1007/978-3-7643-8621-4_15}}.

\bibitem{Whitney1957}
{\sc Whitney, H.}
\newblock {\em Geometric Integration Theory}.
\newblock Princeton University Press, Princeton, N. J., 1957.

\end{thebibliography}

\end{document}